\documentclass[12pt]{article}
\usepackage{amsfonts,amssymb,amsmath,eucal,amsthm,pstricks,pst-plot,pst-node,eucal}
\usepackage[dvips]{graphics}

\newpsobject{showgrid}{psgrid}{subgriddiv=1,griddots=5,gridlabels=0pt}

\newcommand{\cs}{{\varsigma}}
\newcommand{\bA}{\mathsf{A}}
\newcommand{\bT}{\mathsf{T}}
\newcommand{\bE}{\mathsf{E}}
\newcommand{\bM}{\mathsf{M}}
\newcommand{\cE}{\mathcal{E}}
\newcommand{\cS}{\mathcal{S}}
\newcommand{\cI}{\mathcal{I}}
\newcommand{\io}{\vartriangleleft}
\newcommand{\ioe}{\trianglelefteq}
\newcommand{\bo}{\sigma}

\newcommand{\Z}{{\mathbb{Z}}}
\newcommand{\lang}{\left\langle}
\newcommand{\rang}{\right\rangle}

\newcommand{\bG}{\mathsf{G}}
\newcommand{\bH}{\mathsf{H}}
\newcommand{\al}{\alpha} 

\newcommand{\zz}{{\mathfrak{z}}}
\newcommand{\ww}{{\mathbf{w}}}
\newcommand{\ru}{^{\!\!\!\sim}}
\newcommand{\vac}{v_\emptyset}
\newcommand{\cP}{\mathcal{P}}

\newcommand{\bpf}{\begin{proof}}
\newcommand{\epf}{\end{proof}}

\newcommand{\proj}{\mathbf P}

\newcommand{\rarr}{\rightarrow}

\newcommand{\com}{\mathbb{C}}

\newcommand{\pp}{{\mathbf{p}}}

\newcommand{\ul}{\underline}
\newcommand{\sh}{{\textstyle \frac12}}
\newcommand{\nr}[1]{:\!#1\!:}
\newcommand{\gli}{\mathfrak{gl}(\infty)}
\newcommand{\lam}{\lambda}
\newcommand{\la}{\lambda}
\newcommand{\w}{v}

\newcommand{\CS}{\complement S}

\newcommand{\LV}{\Lambda^{\frac\infty2}V}
\newcommand{\LVc}{\Lambda^{\frac\infty2}_0V}

\DeclareMathOperator{\Aut}{Aut}
\DeclareMathOperator{\val}{val}
\DeclareMathOperator{\vrt}{vertices}
\DeclareMathOperator{\Contr}{Contr}
\newcommand{\timeprod}{\operatornamewithlimits{\overrightarrow{\prod}}}

\newtheorem{pr}{Proposition}[section]
\newtheorem{lm}[pr]{Lemma}
\newtheorem{tm}{Theorem}
\newtheorem{cor}[pr]{Corollary}

\theoremstyle{definition}

\numberwithin{equation}{section}

\begin{document}

\title{Virasoro constraints for target curves}
\author{A.~Okounkov and R.~Pandharipande}
\date{August 2003}
\maketitle

\begin{abstract}
We prove generalized Virasoro constraints for the relative
Gromov-Witten theories of all nonsingular target curves.
Descendents of the even cohomology
classes are studied first by localization, degeneration,
and completed cycle methods. Descendents of the
odd cohomology are then controlled by monodromy and geometric
vanishing relations. As an outcome of our results, the relative
theories of target curves are completely and explicitly determined.
\end{abstract}

\setcounter{tocdepth}{2}
\tableofcontents

\setcounter{section}{-1}
\section{Introduction}
\subsection{Overview}

We present here the last in a sequence three papers
devoted to the Gromov-Witten theory of
nonsingular target curves $X$. 
In the
first paper \cite{OP2}, we considered the stationary sector of the
theory formed by the
descendents of the Poincar\'e dual of the
point class. 
The stationary sector was identified in \cite{OP2}
with the Hurwitz theory of $X$ with completed cycles insertions.
In the second paper \cite{OP3}, we found an explicit operator formalism
for the equivariant Gromov-Witten theory of $\proj^1$ in terms
of the  
infinite wedge representation. As a consequence, we proved the
equivariant theory is governed by a 2--Toda hierarchy.

We study here the conjectured Virasoro constraints
for target curves $X$. The standard Virasoro constraints apply 
only to the absolute
Gromov-Witten theory of $X$ and only provide
rules for removing the descendents of the identity class \cite{EgHX}.
The standard 
constraints
are strengthened here in two directions:
\begin{enumerate}
\item[(i)] 
Virasoro constraints 
for
the relative Gromov-Witten theory of target curves $X$ are found,
\item[(ii)]
new constraints providing rules for removing the descendents of
the odd cohomology of $X$ in the relative theory are found.
\end{enumerate}
Our main result is a proof of the strengthened constraints (i)-(ii)
for the relative Gromov-Witten theory of target curves. The
Virasoro conjecture for curves is obtained as a special case of (i).

A complete description of the relative theory of curves is obtained
from the strengthened constraints and the GW/Hurwitz
correspondence of \cite{OP2}. Our main goal in the Introduction
is to present our view of the relative theory of $X$.

\subsection{The relative Gromov-Witten theory of curves}

\subsubsection{}

Let $X$ be a nonsingular target curve 
 of genus $g$. All curves in the paper are projective over $\com$.
Let
\begin{gather}
  1\\
\alpha_1, \ldots, \alpha_g,\ \beta_1, \ldots, \beta_g \\
\omega
\end{gather}
be a basis of
$H^*(X,\com)$ with the following properties:
\begin{enumerate}
\item[(i)] the class $1\in H^0(X,\com)$ is the identity,
\item[(ii)] the classes $\alpha_i\in H^{1,0}(X,\com)$ and
 $\beta_j \in H^{0,1}(X,\com)$
determine a symplectic basis of $H^1(X,\com)$,
$$\int_X \alpha_i \cup \beta_j = \delta_{ij},$$
\item[(iii)] the class $\omega\in H^2(X,\com)$ is the Poincar\'e dual of the point.
\end{enumerate}

\subsubsection{}

We will study the Gromov-Witten theory of $X$ relative to
a finite set of distinct points
$q_1, \ldots, q_{m} \in X$.
Let $\eta^1, \ldots, \eta^m$ be partitions of $d$. The moduli space
$$\overline{M}_{g,n}(X, \eta^1, \ldots, \eta^m)$$ parameterizes
connected, genus $g$, $n$-pointed  stable relative maps with monodromy
$\eta^i$ at $q_i$. Foundational
developments of relative Gromov-Witten theory in
symplectic and algebraic geometry can be found in
\cite{EGH,IP,LR,L}. The absolute Gromov-Witten theory
of $X$ is recovered if $m=0$.

The (nonequivariant) connected Gromov-Witten invariants of $X$ 
relative to $q_1, \ldots, q_m$ are:
\begin{equation}
\label{neddde}
\lang  \prod_{i=1}^{n}
\tau_{k_i}(\gamma_i),\eta^1, \ldots, \eta^m  \rang_{g,d}^{\circ X} =
\int_{[\overline{M}_{g,n}(X,\eta^1,\ldots,\eta^m)]^{vir}}
\prod_{i=1}^n \psi_{i}^{k_i} \,
\text{ev}_i^*(\gamma_i).
\end{equation}
Here,  $\tau_{k}(\gamma)$ denotes the $k$th descendent of the
cohomology class $\gamma\in H^*(X,\com)$. The order of
the descendent insertions in \eqref{neddde} is important for
the classes $\gamma_i$ of odd degree.

The superscript $\circ$ denotes the connected theory.
The corresponding disconnected theory will be denoted by the bracket
$\lang \,\, \rang^\bullet$.  
As we will also use the bare
bracket $\lang\,\, \rang$ for the
disconnected theory, the superscript $\bullet$ will
be used only for emphasis.
We will be primarily interested in the
disconnected theory.

If the target $X$ is $\proj^1$, we will omit the superscript
$\proj^1$.  Most of the  paper will be devoted to the study of  
the relative theories of $\proj^1$ and the elliptic curve $E$.

The subscripted genus may be omitted
in the notation \eqref{neddde} by the dimension constraint in the
nonequivariant theory. If the set of relative points is nonempty,
the subscripted degree may be also omitted.

\subsubsection{}
\label{desp}
We first review the formula for the descendents of $\omega$ 
obtained from the
GW/H correspondence and the theory of completed cycles \cite{OP2}.

Let $d$ be a non-negative integer. Let $\lambda$ be a partition of $d$,
$$\lambda_1 \geq \lambda_2 \geq \lambda_3 \geq \ldots \, .$$
Define the completed cycle $\pp_l(\lambda)$ by the formula:
\begin{equation}\label{defp}
\pp_l(\lambda)=\sum_{i=1}^\infty 
\left[(\lambda_i - i + \tfrac12)^l  - (- i + \tfrac12)^l\right]   + l! c_{l+1}\, ,
\end{equation}
for $l>0$.
The constants $c_{l+1}$ are defined by
$$\sum_{j=0}^\infty c_j z^j = \frac{z/2}{\sinh(z/2)}.$$ 
Notice $\pp_l(\lambda)$ is actually defined by a finite sum for each partition $\lambda$.

Let $\eta$ be a partition of $d$. Let $C_\eta\subset S(d)$ be the
conjugacy class of size $|C_\eta|$ corresponding
to $\eta$. Define the function ${\mathbf f}_\eta(\lambda)$ by
\begin{equation}\label{defff1}
{\mathbf f}_\eta(\lambda)=|C_\eta| \, \frac{\chi^\lambda_\eta}{\dim\lambda},
\end{equation}
where $\chi^\lambda_\eta$ is the character of any element of 
$C_\eta$ in the representation of $S(d)$ corresponding to $\lambda$.

We have the following formula for the descendents of $\omega$ from the GW/H correspondence \cite{OP2}.
\begin{tm}
\label{gwh}
\begin{multline}
  \lang \tau_{z_1}(\omega) \dots \tau_{z_l}(\omega),
  \eta^1,\ldots,\eta^m \rang_{d}^{\bullet X} 
\\= \sum_{|\lambda| =d}
  \left( \frac{\dim\lambda}{d!} \right)^{2-2g} \prod_{i=1}^l
  \frac{\pp_{z_i+1}(\lambda)}{(z_i+1)!}  \prod_{j=1}^m {\mathbf
    f}_{\eta^j}(\lambda) \,.
\end{multline}
\end{tm}

\subsubsection{}

Next, we present our formula governing odd descendents in the
presence of descendents of $\omega$. 

Let $[2k]$ denote the ordered set of integers $(1, \dots, 2k)$.
Let $I(2k)$ denote the set of fixed point free involutions of
$[2k]$. For each element $\sigma \in I(2k)$,
let $o^\sigma_1, \dots, o^\sigma_k$ denote the orbits of $\sigma$.
Each orbit $o^\sigma_i$ is  a two element ordered set $(e_{i1}, e_{i2} )$.
A canonical sign $\epsilon(\sigma)$ is associated to $\sigma$ 
by the parity of the permutation
$$(e_{11},e_{12}, e_{21}, e_{22}, \dots, e_{k1},e_{k2})$$
in the symmetric group $S_{2k}$.

The Gromov-Witten invariants \eqref{neddde} certainly vanish
if an odd number of odd descendent classes are inserted.
Let $\gamma_1, \ldots,\gamma_{2k} \in H^1(X,\com)$ be an even number of
odd classes. 
\begin{tm}
\label{oddd}
\begin{multline}
\lang 
 \prod_{i=1}^{2k}\tau_{y_i}(\gamma_i) 
\prod_{j=1}^l \tau_{z_j}(\omega), \eta^1, \ldots, \eta^m
\rang_{d}^{\bullet X} \\
= \sum_{\sigma \in I(2k)}
\epsilon(\sigma)\  
\prod_{i=1}^k \binom{y_{e_{i1}}+ y_{e_{i2}}}{y_{e_{i1}}}\ \int_X \gamma_{e_{i1}} \cup
{\gamma}_{e_{i2}} \quad  \times \\
\lang 
\prod_{i=1}^k \tau_{y_{e_{i1}}+y_{e_{i2}}-1}(\omega) 
 \prod_{j=1}^l \tau_{z_j}(\omega), \eta^1, \ldots, \eta^m 
\rang_d^{\bullet X} \,.
\end{multline}
\end{tm} 

We easily see the formula is skew-symmetric in the insertions
$\tau_{y_i}(\gamma_i)$. 
Theorem \ref{oddd} will be proven in Section \ref{three} of the paper.

\subsubsection{}
\label{virrr}
Finally, we present our Virasoro constraints for the
Gromov-Witten theory of a genus $g$ curve
$X$ relative to $q_1, \ldots, q_m\in X$.

Let $X^*$ denote the punctured manifold,
$$X^*= X \setminus \{ q_1, \ldots, q_m\},$$
with topological Euler characteristic $\chi(X^*)=2-2g-m$.

We introduce four sets of variables corresponding to the descendents of the
classes 
$1, \alpha_i, \beta_j,$ and $\omega$ respectively:
\begin{gather}
t^0_0, t^0_1, t^0_2,\ldots, \notag \\
s^{i}_0, s^{i}_1, s^{i}_2, \ldots, \ \bar{s}^{j}_0, \bar{s}^{j}_1, \bar{s}^{j}_2, \ldots, 
\label{varrr} \\
t^1_0, t^1_1, t^1_2, \ldots  \notag \, .
\end{gather}
The odd variables $s^i_k, \bar{s}^j_l$ supercommute.
Let $\xi$ denote the formal sum,
$$\xi=\sum_{k\geq 0} t^0_k \tau_k(1) +
\sum_{i=1}^g\sum_{k\geq 0} \left(s^i_k \tau_k(\alpha_i)+ \bar{s}^i_k 
\tau_k(\beta_i)
\right) +\sum_{k\geq 0} t^1_k \tau_k(\omega).$$
Let $Z_d[\eta^1, \ldots,\eta^m]$ be the generating series of disconnected
invariants with fixed relative conditions:
$$Z_d[\eta^1,\ldots, \eta^m] =
\sum_{n \geq 0} \frac{1}{n!}\lang \xi^n, \eta^1, \ldots, \eta^m \rang_d^{\bullet X}.$$
The bracket on the right is expanded multilinearly in the
variables \eqref{varrr}. The supercommutativity of the
odd variables must not be forgotten in the expansion of $Z_d[\eta^1, \ldots, 
\eta^m]$.

We will consider differential operators $D$ in the variables \eqref{varrr} acting
on the series $Z_d[\eta^1, \ldots, \eta^m]$.
The operators will contain only first order derivatives in
the odd variables. The derivative
$$\frac{\partial}{\partial s^i_k} f$$
is defined for monomials $f$ by supercommuting the variable
$s^i_k$ to the left and removing $s^i_k$. The same convention
holds for $\bar{s}^j_l$.

We define the first two Virasoro operators for the
relative theory of $X$ by the following equations:
\begin{eqnarray*}
L_{-1} & = &
-\frac{\partial}{\partial t^0_0} \\
& & +
\sum_{l\geq 0} \Big(t^0_{l+1} \frac{\partial}{\partial t^0_l}
+\sum_{i=1}^g \left(s^i_{l+1} \frac{\partial}{\partial s^i_l}
+\bar{s}^i_{l+1} \frac{\partial}{\partial \bar{s}^i_l}\right)
+ t^1_{l+1} \frac{\partial}{\partial t^1_l}
 \Big) \\ 
& & +
 t^0_0 t^1_0 +  \sum_i s^i_0 \bar{s}^i_0 , 
\end{eqnarray*}

\begin{eqnarray*}
L_0 & = &  -  \frac{\partial}{\partial t^0_1} \\
& & 
 -\chi(X^*) \frac{\partial}{\partial t^1_0} \\
& & +
\sum_{l\geq 0}\Big(  l t^0_l \frac{\partial}{\partial t^0_l}
+
\sum_{i=1}^g \left( (l+1) s^i_l \frac{\partial}{\partial s^i_l}
+  l \bar{s}^i_l \frac{\partial}{\partial \bar{s}^i_l}\right)
+ (l+1)t_l^1 \frac{\partial}{\partial t^1_l} \Big) \\
& &+  \chi(X^*) \sum_{l\geq 0} t^0_{l+1} \frac{\partial}{\partial t^1_l} \\
& & + \frac{\chi(X^*)}{2} t^0_0 t^0_0.
\end{eqnarray*}

The Virasoro operators $L_{-1}$ and $L_0$ for the
relative theory specialize to the corresponding
Virasoro operators for the absolute theory if $m=0$. 
The string, dilaton, and divisor equations for the
relative theory imply the first two Virasoro constraints:
\begin{eqnarray*}
L_{-1} \ Z_d[\eta^1, \ldots, \eta^m] & =& 0, \\
L_0 \ Z_d[\eta^1, \ldots, \eta^m] & = & 0.
\end{eqnarray*}
The derivation of the above constraints is identical
to the corresponding derivation for the absolute theory.

In the definition of the remaining 
 Virasoro operators, we will 
 use the Pochhammer symbol,
$$(a)_{b} = \frac{(a+b-1)!}{(a-1)!}.$$
For $k>0$, the operators $L_k$ for
the
relative theory of $X$ are defined by:
\begin{eqnarray*}
L_k &  = & - (k+1)! \ \frac{\partial}{\partial t^0_{k+1}} \\
&  & - \chi(X^*)\ (1)_{k+1}\ \sum_{r=1}^{k+1}{\frac{1}{r}}\ 
\frac{\partial}{\partial t^1_k} \\
& & + \sum_{l\geq 0} \left(
   (l)_{k+1}\
t^0_l \frac{\partial}{\partial t^0_{k+l}} +
(l+1)_{k+1}\
t^1_l \frac{\partial}{\partial t^l_{k+l}} \right) \\
& & + \sum_{l\geq 0}\sum_{i=1}^g \left(
(l+1)_{k+1}\
s^i_l \frac{\partial}{\partial s^i_{k+l}} +
(l)_{k+1} \
\bar{s}^i_l \frac{\partial}{\partial \bar{s}^i_{k+l}}\right)  \\
& & + \chi(X^*) \sum_{l\geq 0} 
 (l)_{k+1} \ \sum_{r=l}^{k+l} \frac{1}{r} \ 
t^0_{l} \frac{\partial}{\partial t^1_{k+l-1}} \\
& & +  \frac{\chi(X^*)}{2}\sum_{l \geq 0}^{k-2} 
(l+1)! (k-l-1)!  \frac{\partial}{\partial t^1_l}
\frac{\partial}{\partial t^1_{k-l-2}}.
\end{eqnarray*}

The operators $L_k$ are easily seen to satisfy the Virasoro bracket,
$$[L_n,L_m] = (n-m) L_{n+m},$$
and thus determine a representation of the subalgebra of
the Virasoro algebra spanned by  
holomorphic vector fields,
$${\mathcal V}=
\left\{ -z^{k+1}\frac{\partial}{\partial z}\right\}_{k\geq -1}.$$

A central result of the paper is a proof of the Virasoro
constraints for the relative theory of $X$.
\begin{tm} 
\label{rv}
For all $k\geq -1$,
$\ L_k  Z_d[\eta^1, \ldots, \eta^m]=0.$
\end{tm}
The proof is presented in two parts. The Virasoro 
constraints for the descendents of the even cohomology
classes of $X$ are proven first in Sections \ref{two}- \ref{llee}. The
full constraints are established in Section \ref{three}.
The Virasoro constraints for the absolute
theory of $X$ are obtained if $m=0$.

Theorems \ref{gwh} -- \ref{rv}
uniquely determine the
relative Gromov-Witten theory of $X$. 
Every relative invariant of $X$ can be
efficiently calculated.

\subsubsection{}
The Virasoro constraints for the relative Gromov-Witten theory
of $X$ represent a strengthening of the standard Virasoro constraints.
The operator $L_k$ provides a rule for the removal of the descendent
$\tau_k(1)$ in the relative theory of $X$.

We define additional differential operators $D^i_k$ and $\bar{D}^i_k$
for $k\geq -1$ by:
\begin{align*}
  D^i_k = -(k+1)! \, \frac{\partial}{\partial s^i_{k+1}} + 
& \sum_{l=0}^\infty (l)_{k+1} \, t^0_l \, \frac{\partial}{\partial s_{k+l}}\\
+ &
\sum_{l=0}^\infty (l+1)_{k+1} \, \bar {s}^i_l \, \frac{\partial}{\partial t^1_{k+l}}
\end{align*}
\begin{align*}
  \bar{D}^i_k = -(k+1)! \, \frac{\partial}{\partial \bar{s}^i_{k+1}} 
+ & 
\sum_{l=0}^\infty (l)_{k+1} \, t^0_l \, \frac{\partial}{\partial \bar{s}^i_{k+l}}\\
-&
\sum_{l=0}^\infty (l+1)_{k+1} \,  s^i_l \, \frac{\partial}{\partial t^1_{k+l}}
\end{align*}
These operators annihilate the generating series $Z_d[\eta^1,\ldots, \eta^m]$
and provide rules for the removal of the descendents $\tau_k(\alpha_i)$
and
$\tau_k(\beta_i)$.

\begin{tm} For all $k\geq -1$,
\label{oddop}
\begin{eqnarray*}
D^i_{k} \ Z_d[\eta^1, \ldots, \eta^m] & =& 0, \\
\bar{D}_k^i \ Z_d[\eta^1, \ldots, \eta^m] & = & 0.
\end{eqnarray*}
\end{tm}

Theorem \ref{oddop} is derived from Theorems \ref{gwh} -- \ref{rv} in
Section \ref{three} and represents our second strengthening of the
standard Virasoro constraints.

\subsubsection{}

The operators $L_k, D^i_k, \bar{D}^i_k$ satisfy the following 
commutation relations:
\begin{align*}
 [L_n, L_m] & =  (n-m) L_{n+m},\\
           [L_n, D^i_m] & =  -(m+1) D^i_{n+m},\\
 [L_n, \bar{D}^i_m] & =  (n-m) \bar{D}^i_{n+m}.\\
\end{align*}
The odd operators all anti-commute:
$$\{D^i_n, D^j_m\}= \{D^i_n, \bar{D}^j_m\}= \{\bar{D}^i_n,
\bar{D}^j_m\}=0.$$

Let the  operators $\{L_k\}_{k\geq -1}$ be identified with the
Lie algebra of holomorphic vector field $\mathcal V$.
Then, the  operators $\{ D^i_{k} \}_{k\geq-1}$ define
a $\mathcal V$-module isomorphic to 
$$\{ -z^{k+1} \}_{k \geq -1}$$
with the action defined by differentiation, and
the operators $\{\bar{D}^i_k\}_{k\geq -1}$ define
a $\mathcal V$-module isomorphic to the
adjoint representation.

\subsection{Plan of the paper}
The Virasoro constraints for the even theory are studied by degeneration
in Section \ref{two}.
The basic building blocks of the degeneration scheme 
are the cap, the tube, and the pair of pants.
From the algebraic perspective, the building blocks
may be viewed as $\proj^1$ relative to 1, 2, and 3 points respectively.
The main result of Section \ref{two} is
the
reduction of the even Virasoro conjecture to the case of the cap.

The theory of the cap is studied by localization in
Section \ref{one}. 
The theory of the tube arises in the vertex
integrals of the localization formula for the cap. 
The relationship between the theories of the cap and tube
plays an essential role in our proof of the even Virasoro
constraints.

In Section \ref{oopp}, the theory of the cap is
expressed in terms of vacuum expectations of  operators  in
the infinite wedge representation $\LV$.
The Virasoro constraints for the cap are derived
in Section \ref{llee} using the operator formalism.

The  
full relative theory including the odd classes
is studied in Section \ref{threeodd}.
Several techniques including monodromy invariance
and geometric vanishings are used. The proofs
of Theorems \ref{oddd} -- \ref{oddop} are
completed in Section \ref{three}.

\subsection{Acknowledgments}
We thank J.~Bryan, E.~Getzler, and 
T.~Graber for discussions about
Gromov-Witten theory and the Virasoro constraints.

A.~O.\ was partially supported by 
DMS-0096246 and fellowships from the Sloan and Packard foundations.
R.~P.\ was partially supported by 
DMS-0236984 and fellowships from the Sloan and Packard foundations.
The research presented here was partially pursued during a visit by R.~P. to
MSRI in the spring of 2002.   
The paper was written in the summer of 2003 in Princeton and Lisbon.

\pagebreak

%
%

\section{Virasoro constraints
for even classes}

\label{two}

\subsection{Overview}

In Sections \ref{two} -- \ref{llee},
we will consider  the relative Gromov-Witten theory of
$X$ with only the descendents $\tau_k(\gamma)$ of the even
cohomology classes $$\gamma\in H^{2\bullet}(X,\com).$$ 
The odd theory will be studied in
Sections \ref{threeodd} and \ref{three}.

The Virasoro constraints for the relative Gromov-Witten
theory of $X$
are easily seen to respect the even classes. 
Let $\xi^{*}$ denote the even sum,
$$\xi^*=\sum_{k\geq 0} t^0_k \tau_k(1) +
\sum_{k\geq 0} t^1_k \tau_k(\omega),$$
and let
Let $Z^*_d[\eta^1, \ldots,\eta^m]$ be the generating series of 
even relative invariants:
$$Z^{*}_d[\eta^1,\ldots, \eta^m] =
\sum_{n \geq 0} \frac{1}{n!}\lang (\xi^{*})^n, 
\eta^1, \ldots, \eta^m \rang_d^{\bullet X}.$$
Let $L^*_k$ denote the restricted Virasoro 
operators,
$$L^*_k = L_k |_{\{s^i_p, \bar{s}^i_q =0\}}.$$
The full Virasoro constraints {\em imply} the
even Virasoro constraints:
$$L_k^* \ Z^*_d[\eta^1, \ldots, \eta^m]=0,$$
for all $k\geq -1$.
Our goal in Sections \ref{two} -- \ref{llee}
is to prove the even Virasoro constraints
for all relative target curves $X$.

\subsection{Virasoro reactions}
\subsubsection{}
The Virasoro constraints provide rules for removing
$\tau_k(1)$ insertions in the relative Gromov-Witten
theory of target curves $X$. We will describe the Virasoro
rule for the removal of $\tau_k(1)$ from
\begin{equation}
\label{prere}
\lang \tau_k(1) \prod_{i} \tau_{l_i}(\gamma_i) \rang^X
\end{equation}
as a {\em reaction}.

The descendent $\tau_k(1)$ is viewed as unstable and
subject to decay. The descendent $\tau_k(1)$
decays via five type of reactions with the other
insertions of \eqref{prere}. These reactions are
detailed in the table below. The columns show the number,
the type, and the actual formulas for the
reactions. 

\renewcommand{\arraystretch}{2}
\begin{center}
\begin{tabular}{|l|l|rl|}
\hline 
(i) & $1+1\to1$ & 
$\displaystyle\tau_k(1) \, \tau_l(1) \to$&$\displaystyle 
\binom{k+l-1}{k} \, \tau_{k+l-1}(1)$ \\
(ii) & $1+1\to\omega$ & 
$\displaystyle\tau_k(1) \, \tau_l(1) \to$&$\displaystyle 
\chi(X^*)\,\binom{k+l-1}{k} \, \left(\sum_{j=l}^{k+l-1}\frac{1}{j}\right)
\,\tau_{k+l-2}(\omega)$ \\
(iii) & $1+\omega\to\omega$ & 
$\displaystyle\tau_k(1) \, \tau_l(\omega) \to$&$\displaystyle 
\,\binom{k+l}{k} \, 
\,\tau_{k+l-1}(\omega)$ \\
(iv) & $1\to\omega$ &  
$\displaystyle\tau_k(1) \to$&$\displaystyle 
-\chi(X^*)\,\left(\sum_{j=1}^{k}\frac{1}{j}\right)
\,\tau_{k-1}(\omega)$ \\
(v) & $1\to\omega+\omega$ &  
$\displaystyle\tau_k(1) \to$&$\displaystyle 
\frac{\chi(X^*)}{2k}\,\sum_{i=0}^{k-3} \binom{k-1}{i+1}^{-1} 
\tau_{i}(\omega) \, \tau_{k-i-3}(\omega)$
\\
\hline
\end{tabular} 
\end{center}

Recall $X^*$ is the manifold obtained by removing the
relative points of $X$, and $\chi(X^*)$ is the topological 
Euler characteristic.
The reactions (ii), (iv), (v), whose intensity 
involves a factor of $\chi(X^*)$, are \emph{extensive}
reactions. 
The remaining reactions (i) and (iii) 
are  \emph{intensive} reactions.

The Virasoro rule for the removal of $\tau_k(1)$ 
from \eqref{prere} for $k\geq 1$ is to sum over all the
invariants arising as outputs of the the five decay
reactions. For example,
$$\lang \tau_2(1) \tau_3(\omega)\rang^X   =   \binom{5}{2}\lang
 \tau_4(\omega) \rang^X
 - \frac{3}{2}\, \chi(X^*) \lang \tau_1(\omega) \tau_3(\omega)\rang^X,$$
by rules (iii) and (iv).

The five reactions cover all the terms of the
even Virasoro operators with exception of the constant terms
in $L^*_{-1}$ and $L^*_0$.
The Virasoro rules are therefore equivalent to the
even Virasoro constraints for $k\geq 1$.
We will prove the Virasoro rules are valid for the
even relative Gromov-Witten theory of $X$.

For $\tau_0(1)$ and $\tau_1(1)$ the Virasoro rules
can be supplemented to incorporated the constant terms.
However, since the Virasoro constraints  $L^*_{-1}$ and $L^*_{0}$
are proven, we will not investigate them further.

\subsubsection{}

The stationary theory of $X$ relative to $q_1, \ldots, q_m$
is determined by Theorem \ref{gwh} via the GW/H correspondence.
The Virasoro constraints uniquely determine an 
even theory from the stationary sector.

\begin{pr}
\label{hhh}
There exists a unique solution to the even Virasoro constraints,
$$ L_k^* \bar{Z}_d^*[\eta^1,\ldots,\eta^m] =0,\ \  \forall k\geq -1,$$ 
which extends the stationary Gromov-Witten theory of $X$.
\end{pr}

\bpf
The coefficients of $\bar{Z}_d^*[\eta^1,\ldots,\eta^m]$ determine
a new bracket
$$\lang \prod_{i} \tau_{k_i} (\gamma_i) \rang ^{-},$$
for even classes $\gamma_i$.
The solution is said to {\em extend} the stationary theory of $X$ { if} the bracket 
$\lang, \rang^-$
agrees with the relative Gromov-Witten
bracket $\lang, \rang$ in case all insertions
are descendents of $\omega$.

The uniqueness of the  solution $\tilde{Z}_d^*[\eta^1,\ldots,\eta^m]$
is clear from the Virasoro reactions. After repeated applications, the
reactions remove all the
descendents of the identity class from $\lang, \rang^-$ and leave only
the stationary descendents.

To prove existence, we must prove the Virasoro rules are
compatible. Given a bracket 
 $$\lang \prod_{i} \tau_{k_i} (\gamma_i) \rang ^{-},$$
we must prove the reduction of the bracket to the
stationary theory is independent of the {\em order} of application
of the Virasoro rules.

The compatibility is easily obtained by induction on
the number of 
descendents of the identity in the bracket $\lang, \rang^{-}$
and the commutation relation 
$$[L^*_n,L^*_m]=(n-m) L^*_{n+m}$$
of the Virasoro operators.
\epf

\subsection{Degeneration}

\subsubsection{}
Let $X$ be a target curve with relative points $q_1, \ldots, q_m$.
We will consider nodal degenerations of $X$ of two types:
\begin{enumerate}
\item[(i)] $X$ degenerates to $X' \cup X''$ intersecting in a
           node $q_*$. The relative points are distributed  in the degeneration
           to $q_1',\ldots,q_{m'}'$ on $X'$ and
           $q''_{1}, \ldots, q''_{m''}$ on $X''$.
\begin{figure}[!htbp]
  \begin{center}
    \scalebox{0.64}{\includegraphics{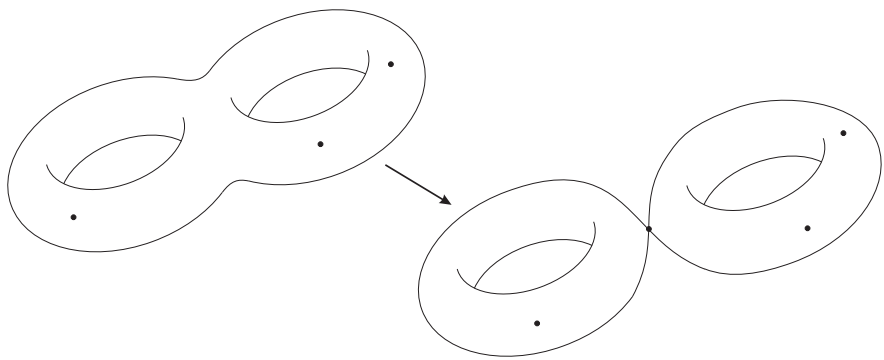}} 
    \caption{Nodal degeneration of type (i)}
    \label{f1}
  \end{center}
\end{figure}
\item[(ii)] $X$ degenerates to an irreducible curve $X'$ of 
            geometric genus 
$$g(X')= g(X)-1$$
            with node $q_*$.
\begin{figure}[!htbp]
  \begin{center}
    \scalebox{0.64}{\includegraphics{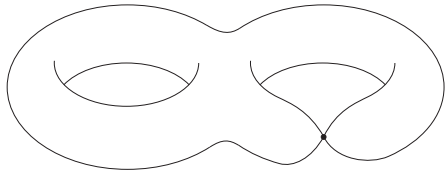}} 
    \caption{Nodal degeneration of type (ii)}
    \label{f2}
  \end{center}
\end{figure}

\end{enumerate}

The degeneration formula expresses the relative Gromov-Witten invariants
of $X$ in terms of the relative theory of the degenerations.
Consider the relative invariant
\begin{equation}
\label{tt1}
\lang \prod_{i=1}^k \tau_{y_i}(1) \prod_{j}^l \tau_{z_j}(\omega),
\eta^1, \ldots, \eta^m  \rang^X. 
\end{equation}
Let $T\subset \{1, \ldots, l\}$ be a subset.
The type (i) degeneration formula for the invariant \eqref{tt1} is:
\begin{multline}
\label{degen1} 
\sum_{S\subset \{1,\ldots,k\}} \sum_{|\mu|=d}\ \  
\lang \prod_{i\in S} \tau_{k_i}(1) \prod_{j\in T} \tau_{l_j}(\omega),
\eta^1, \ldots, \eta^{m'}, \mu  \rang^{X'} 
\zz(\mu) \quad \times \\
 \lang \prod_{i\notin S} \tau_{k_i}(1)
 \prod_{j\notin T} \tau_{l_j}(\omega),
\eta^1, \ldots, \eta^{m''}, \mu  \rang^{X''}\,. 
\end{multline}
The type (ii) degeneration formula for the invariant \eqref{tt1} is:
\begin{equation}
  \label{degen2}
  \sum_{|\mu|=d} \zz(\mu)
\lang \prod_{i} \tau_{k_i}(1) \prod_{j} \tau_{l_j}(\omega),
\eta^1, \ldots, \eta^n, \mu, \mu  \rang^{X'} \,. 
\end{equation}
Here, the 
automorphism factor $\zz(\mu)$ is defined by:
\begin{equation}
  \label{zzmu}
  \zz(\mu)= \left|\Aut(\mu)\right| \prod_{i=1}^{\ell(\mu)} \mu_i\,,
\end{equation}
where $\Aut(\mu)$
is the symmetry group permuting equal parts of $\mu$. 

Proofs of the degeneration formulas \eqref{degen1} and 
\eqref{degen2}  can be found in \cite{IP,LR,L}.

\begin{pr}
\label{zzz}
If the relative theories of the degenerations of either type
satisfy the even Virasoro constraints, then the original relative
theory of $X$ satisfies the even Virasoro constraints.
\end{pr}

\bpf
The Proposition is obtained by an elementary verification
of the compatibility of the Virasoro rules with the
degeneration formulas.
\epf

\subsubsection{}
The cap $\mathbf C$ is the target determined $\proj^1$ relative to $\infty\in \proj^1$.
We will see the relative theory of the cap governs the even 
relative theories all of target curves.
Motivated by the operator formalism of the infinite wedge
representation,
we will denote the relative invariants of $\mathbf C$ by:
\begin{equation*}
\lang  \left. \prod_{i=1}^{n}
\tau_{k_i}(\gamma_i) \right|  \eta  \rang .
\end{equation*}
The bracket $\lang\,\, \rang^{\mathbf C}$ may be used to emphasize the target.
\subsubsection{}

Let $X$ be a target curve with relative points $q_1, \ldots,q_m$.
Consider the type (i) degeneration of $X$ to $X \cup {\mathbf C}$,
\begin{figure}[!htbp]
  \begin{center}
    \scalebox{0.64}{\includegraphics{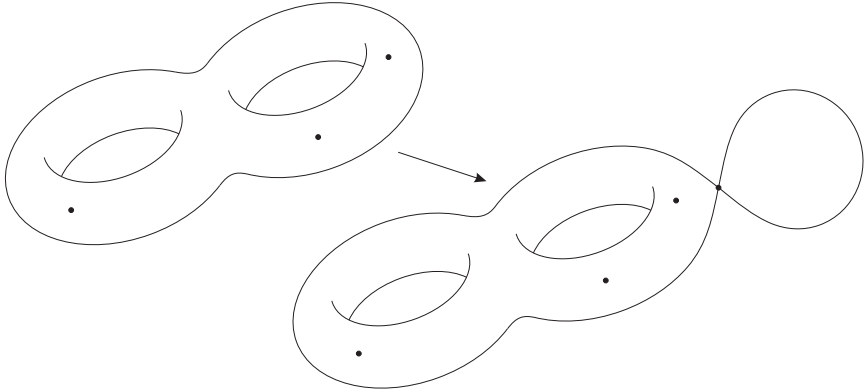}} 
    \caption{Degeneration used in the proof of Proposition \ref{fff}}
    \label{wpwz}
  \end{center}
\end{figure}
where ${\mathbf C}$ is the cap and all the relative points $q_1, \ldots, q_m$ remain on $X$, 
see Fig.~\ref{wpwz}. 
We will use the degeneration formulas and the GW/H correspondence to
prove the following result.

\begin{pr}
\label{fff}
The even theory of $X$ relative to $q_1, \ldots, q_m$ and the
theory of the cap uniquely determine the
even theory of $X$ relative to $q_1, \ldots, q_m,q_*$.
\end{pr}

\bpf
Let $p(d)$ be the number of partitions of size $d$.
Consider the $\infty \times p(d)$ matrix $M_d$, indexed by 
monomials $L$ in the descendents of $\omega$ and partitions $\mu$ of $d$,
 with
coefficient
$$\lang L \ |\  \mu \rang^{\mathbf C}_d$$
in position $(L,\mu)$.
Since the completed cycles, $$\pp_1, \pp_2, \pp_3, \ldots,$$ generate the
algebra of shifted symmetric functions, the matrix
$M_d$ has full rank equal to $p(d)$ by the GW/H correspondence, see \cite{OP2}.  

We prove the Proposition by induction on the number of 
descendents of the identity in the bracket 
$\lang \,\, \rang ^{X}$
relative to $m+1$ points. The base case and the induction step are
proven simultaneously.

Consider the invariants of $X$ relative to $m+1$ points,
\begin{equation}
\label{yone}
\lang  
\prod_i \tau_{k_i}(1) \prod_j\tau_{l_j}(\omega), \eta^1, \ldots, \eta^m, \mu\rang _d^X,
\end{equation}
defined by fixing the descendent insertions and varying  $\mu$
among all partitions of $d$.

We will determine the invariants \eqref{yone} from the
invariants of $X$ relative to $m$ points,
\begin{equation}
\label{yall}
\lang L\ 
\prod_i \tau_{k_i}(1) \prod_j\tau_{l_j}(\omega), \eta^1, \ldots, \eta^m\rang_d^X,
\end{equation}
defined by all monomials $L$  in the descendents of $\omega$.

We apply the  type (i) formula 
for the degeneration of Fig.~\ref{wpwz} to the invariants \eqref{yall}. 
All  the descendents of $\omega$
remain on $X$ in the degeneration except for those 
in $L$ which distribute to ${\mathbf C}$. 
By induction, we need only analyze the terms of the degeneration formula
in which the descendents of the identity remain on $X$. Then,
since $M_d$ has full rank, the invariants \eqref{yone} are determined
by the invariants \eqref{yall}.
\qed 
\vspace{+10pt}

The {\em depth} $r$ theory of $X$ relative to $q_1, \ldots, q_m$
consists of the even invariants with at most $r$ descendents of the
identity class. In particular, the depth $0$ theory coincides with the
stationary theory of $X$.
The proof of Proposition \ref{fff} yields a refined
result: the depth $r$ theory of $X$ relative to $m$ points
and the depth $r$ theory of the cap uniquely determine the
depth $r$ theory of $X$ relative to $m+1$ points.

\subsubsection{}

Proposition \ref{fff}  implies the 
theory of the cap determines the relative theory of all
curves. 

\begin{pr} 
\label{eee}
If the theory of the cap satisfies
the Virasoro constraints, then the relative theories of all curves
satisfy the even Virasoro constraints.
\end{pr}

\bpf
Assume the relative Gromov-Witten theory of the cap satisfies the Virasoro
constraints. 
We will prove by induction on $m$ that the
Gromov-Witten theory of $\proj^1$ relative
to $m$ points also satisfies the
Virasoro constraints. The Proposition is  then a consequence of 
Proposition \ref{zzz} since every curve $X$ admits a sequence of type (ii)
degeneration to $\proj^1$.
The base of the induction is satisfied since $\proj^1$ relative to
1 point is the cap.

Let $m>1$.
The Virasoro constraints define
a unique extension of the stationary
Gromov-Witten theory of $\proj^1$ relative $m$ points. 
By Proposition \ref{zzz} and the induction
hypothesis, the extension defined
by the Virasoro constraints
is compatible with the degeneration formulas.

By Proposition \ref{fff}, the Gromov-Witten theory 
of $\proj^1$ relative $m$ points is uniquely determined
by the Gromov-Witten theory of $\proj^1$ relative $m-1$ points
by considering the type (i) degeneration of $\proj^1$ to 
$\proj^1 \cup {\mathbf C}$. Therefore, the
extension defined by the Virasoro constraints
coincides with the Gromov-Witten theory of $\proj^1$
relative to $m$ points.
\epf

The proof of Proposition \ref{eee} yields a refined
result: the Virasoro constraints
for the depth $r$  theory of the cap implies
the Virasoro constraints for the depth $r$ relative theories of
all target curves.

\pagebreak

\section{The equivariant theory of the  cap}
\label{one}
\subsection{Overview}

The theory of the cap governs  the even relative theory of all target curves. 
We will now study the equivariant 
relative Gromov-Witten theory of the cap
by localization.
The results will be cast in the operator formalism
of $\LV$ in Section \ref{oopp}.

\subsection{Theories of the tube}
The tube ${\mathbf T}$ is the target determined by $\proj^1$ relative
to $0,\infty \in \proj^1$. The automorphism group of the tube fixing the 
relative points is $\com^*$.
There are two relative theories of the
tube:
\begin{enumerate}
\item[(i)] The {\em parameterized} theory 
concerns integration
over the moduli space $\overline{M}_{g,n}(\proj^1, \mu,\nu)$.
\item[(ii)] The {\em unparameterized} theory concerns
integration over the moduli space
 $\overline{M}^\sim_{g,n}(\proj^1, \mu,\nu)$
obtained by identifying maps which differ by an automorphism of the
target.
\end{enumerate}
The relative Gromov-Witten theory of target curves specializes to
the parameterized theory of the tube.

The unparameterized theory is special to the geometry of the
tube.
Let $$U_{g,n}(\proj^1,\mu,\nu)\subset \overline{M}_{g,n}(\proj^1,\mu,\nu)$$ be the
locus of maps with finite $\com^*$-stabilizers for the induced
action on the moduli of maps. 
The moduli space $\overline{M}_{g,n}^\sim(\proj^1,\mu,\nu)$ is defined by:
\begin{equation}
\label{ned}
\overline{M}^\sim_{g,n}(\proj^1,\mu,\nu) = U_{g,n}(\proj^1,\mu,\nu)\ / \ \com^*.
\end{equation}
Since the tube dilates with the $\com^*$-action,
we view $\overline{M}^\sim_{g,n}(\proj^1,\mu,\nu)$ as a moduli space of
maps to {\em rubber}.

The unparameterized theory of the tube arises naturally in the vertex
integrals of the localization formula for the cap.

\subsection{Localization for the cap}
\label{lcap}

\subsubsection{}
Let $V= \com \oplus \com$.
Let the algebraic torus $\com^*$ act on $V$ with weights $(0,1)$:
$$
\xi \cdot(v_1,v_2)= (v_1, \xi \cdot v_2)\,.
$$
Let $\proj^1$ denote the projectivization $\proj(V)$.
There is a canonically induced
$\com^*$-action on $\proj^1$ with fixed points $0,\infty \in \proj^1$.

The $\com^*$-equivariant cohomology ring of a point is $\com[t]$ where
$t$ is the first Chern class of the standard representation.
The $\com^*$-equivariant cohomology ring $H^*_{\com^*}(\proj^1, \com)$
is canonically a $\com[t]$-module.
The (localized) equivariant cohomology of $\proj^1$ is
spanned by the classes of the fixed points,
$$[0],[\infty] \in H^*_{\com^*}(\proj^1, \com).$$

\subsubsection{}
We will study here the equivariant relative Gromov-Witten theory
of the cap.
Let $\infty\in \proj^1$ be the relative point of the cap.
The $\com^*$-action lifts to the moduli space $\overline{M}^\bullet_{g,n}(\proj^1,\nu)$.
As the virtual structure is canonically equivariant,
we may define equivariant relative invariants by equivariant
integration:
$$
\lang \left. \prod_{i=1}^n \tau_{k_i}(0)  \, \prod_{j=1}^m
 \tau_{l_j}(\infty) \right|
\nu \rang_{g} = \int_{[\overline{M}^\bullet_{g,n+m}(\proj^1,\nu)]^{vir}} \prod_{i=1}^n
\psi_i^{k_i} \text{ev}^{*}([0]) 
\prod_{j=1}^m
\psi_j^{l_j} \text{ev}^{*}([\infty]). 
$$
The equivariant relative invariants take values in
$\com[t]$.

The virtual class is defined by an (equivariant)
perfect obstruction theory \cite{GV, L}. The virtual
localization formula of \cite{GP} is applied to the
case of relative maps in \cite{GV}. 
A presentation of the localization formula for
relative maps to $\proj^1$ can also be
found in \cite{FabPan}.

\subsubsection{}

The virtual localization formula for the cap will be used
to study the equivariant relative invariants. 
The formula is best expressed in
terms of the generating series of equivariant
relative invariants, 
\begin{multline} \label{bG}
  \bG(z_1,\dots,z_n,w_1,\dots,w_m |\nu) = \\ 
\sum_{g} \sum_{k_i} \sum_{l_j} \, \prod_{i=1}^n z_i^{k_i+1} 
\prod_{j=1}^m w_j^{l_j+1} \,
\lang \left. \prod \tau_{k_i}(0)  \, \prod \tau_{l_j}(\infty) 
\right| \nu \rang_g .  
\end{multline}
As disconnected invariants are considered, the genus $g$ may be negative
in the outer sum. The genus variable $u$ of \cite{OP3} is omitted
here without any loss of information.

The nonequivariant
specialization of $\bG$ will be denoted by $\bG'$: 
\begin{multline}\label{bG'}
  \bG'(\left. x_1,\dots,x_n,y_1,\dots,y_m \right|\nu) = \\
\sum_g \sum_{k_i\ge 0} \sum_{l_j\ge 0} \, \prod_{i=1}^n x_i^{k_i+1} 
\prod_{j=1}^m y_j^{l_j+1} \, 
\lang \left. \prod \tau_{k_i}(\omega)  \, \prod \tau_{l_j}(1) 
\right| \nu \rang_g  \,.
\end{multline}
By the dimension constraint,
the genus $g$ may be omitted in the above formula.

The function $\bG'$ is obtained from the function $\bG$ by using the
relations
$$
\omega = [0],$$
\begin{equation}
1 = \frac{[0]-[\infty]}{t} \,,
\label{10i}
\end{equation}
and then letting $t\to 0$. The procedure is very 
much like differentiating $\bG$ with 
respect to $t$. 

The depth $r$, nonequivariant, relative theory of the cap is
determined by the set of functions, 
$$\{ \bG(z_1,\ldots,z_n, w_1, \ldots, w_m|\nu)\}_{n< \infty, m \leq r },$$
by nonequivariant specialization.

\subsubsection{}

To achieve uniformity in the localization formulas and the 
operator formalism, we will include 
unstable contributions in the generating functions such as
\eqref{bG}. We follow the conventions of \cite{OP3}
concerning the unstable contributions. The true 
Gromov-Witten invariants are obtained from the coefficients of
terms of positive degree in all the variables.

\subsubsection{}

The localization formula expresses $\bG$ in terms of Hodge integrals 
over 0 and rubber integrals over $\infty$.
We first introduce the generating functions associated to
these vertex integrals.

Let $\bH(z_1, \ldots, z_n,u)$ denote the disconnected 
$n$-point function of $\lambda$-linear Hodge integrals 
defined in \cite{OP3}. We follow here
the Hodge integral conventions of \cite{OP3} governing the unstable
cases. The function $\bH$ is identified as a vacuum expectation in
$\LV$ in \cite{OP3} and is fully determined.

The rubber integrals which arise over $\infty$ are of the
following form:
\begin{equation*}
\lang \mu, k \left|  \prod_{i=1}^n \tau_{k_i}\right|  \nu
\rang^{\sim}_g
= \int_{[\overline{M}^{\bullet\sim}_{g,n}(\proj^1,\mu,\nu)]^{vir}} \psi^k 
\prod_{i=1}^n \psi_i^{k_i},
\end{equation*}
where $\psi$ is the cotangent line to
the target at the relative point.
The superscript $\sim$ indicates the rubber target.
Let
\begin{multline}\label{rG}
  \bG^\sim(\mu,s\ |\ w_1,\dots,w_m\ |\ \nu) = \\
\sum_g \sum_{k\geq 0} \sum_{l_j\ge 0} \, s^{k-m+1} \prod_{j=1}^m w_j^{l_j+1} 
\, 
\lang \mu,k \left| \prod \tau_{l_j}
\right| \nu \rang_g\ru \,.
\end{multline}

A direct application of the virtual localization formula for the
cap yields the following result.

\begin{pr}\label{lcapp}
$\bG(z_1,\dots,z_n,w_1,\dots,w_m |\nu) =$ 
$$
\sum_{|\mu|=d} 
\frac{1}{\zz(\mu)}
\frac{
    (\tfrac 1t)^{\ell(\mu)}} {t^{d+n}}
 \, \prod
    \frac{\mu_i^{\mu_i}}{\mu_i !}  \bH(\mu,tz,\tfrac 1t)  \  
  \zz(\mu) \ \bG^\sim(\mu,-\tfrac 1t\ |\ w_1,\dots,w_m\ | \ \nu)
.
$$
\end{pr}

\noindent Here, $\zz(\eta)$ is the automorphism factor
defined in \eqref{zzmu}. 

The Proposition will be expressed in terms of the operator formalism of
the infinite wedge representation in Section \ref{two}.

\subsection{Rubber calculus I}
\label{pu}

\subsubsection{}
The  rubber integrals,
\begin{equation}
\label{ep1}
\lang \mu, k \left|  \prod_{i=1}^n \tau_{k_i}\right|  \nu
\rang^{\sim}_g
= \int_{[\overline{M}^{\bullet \sim}_{g,n}(\proj^1,\mu,\nu)]^{vir}} \psi^k 
\prod_{i=1}^n \psi_i^{k_i},
\end{equation}
arise in the localization formula for the cap.
The relative Gromov-Witten invariants of the tube are:
\begin{equation}
\label{ttt}
\lang \mu \left|  \prod_{i=1}^n \tau_{k_i}(\gamma_i) \right|  \nu
\rang_g
= \int_{[\overline{M}^\bullet_{g,n}(\proj^1,\mu,\nu)]^{vir}} 
\prod_{i=1}^n \psi_i^{k_i} \text{ev}_i^*(\gamma_i),
\end{equation}
where $\gamma_i \in H^*(\proj^1, \com)$.
It will be important for us to express the rubber integrals \eqref{ep1} in
terms of the tube theory \eqref{ttt} {\em without} cotangent line classes
at the relative points.

If the marking set is nonempty,  the rubber integral \eqref{ep1} is immediately
expressed as a tube integral with cotangent line classes at the relative points.

\begin{lm}
For $n>0$,
\label{sittt}
\begin{equation}
\lang \mu, k \left| \prod_{i=1}^n \tau_{k_i} \right|  \nu
\rang^{\sim} =
\lang \mu,k \left|  \tau_{k_1}(\omega) 
\prod_{i=2}^n \tau_{k_i}(1) \right|   \nu
\rang.
\end{equation}
\end{lm}

\bpf
The conceptually simplest proof of the Lemma is by a direct identification 
of moduli spaces,
$$\overline{M}^{\bullet \sim}_{g,n}(\proj^1, \mu,\nu)
 \stackrel{\sim}{=} 
\text{ev}_1^{-1}(x) \subset \overline{M}^{\bullet}_{g,n}(\proj^1,\mu,\nu),$$
by the quotient description \eqref{ned}.
Here, $x\in \proj^1$ is any point not equal to 0 or $\infty$.
The Lemma then follows after
an identification of obstruction theories and virtual classes.

Alternatively, a
formal derivation of the Lemma from a localization calculation of
the tube integral can be found. The details are left to the
reader.
\epf

The cotangent line $\psi$ at the relative point $q$
is easily analyzed on the moduli
space $\overline{M}_{g,n}^\bullet(\proj^1,\mu,\nu)$.
Since  $q$ lies over the
fixed point $0\in \proj^1$, we obtain
\begin{equation}
\label{sde}
\psi= \psi_0 + \Delta,
\end{equation}
where $\psi_0$ is the trivial cotangent line at $0\in \proj^1$ and
$\Delta$ is the divisor corresponding to proper degenerations
in the Artin stack of degeneration of 
the relative space $0\in \proj^1$. 

Lemma \ref{sittt}, equation \eqref{sde},
 and the splitting formula for the relative theory together
yield
the following recursive relation in case $n>0$ and $k>0$:
\begin{multline}
  \label{rcr}
 \lang \mu,k \left| \prod_{i=1}^n \tau_{k_i} \right|
\nu\rang\ru = \sum_{S\subset\{2,\dots,n\}} \sum_{\eta} \\
\lang \mu,k-1 \left| \prod_{i\notin S} \tau_{k_i} \right|
\eta\rang\ru\, \zz(\eta) \,
\lang \eta \left| \tau_{k_1}(\omega) \, 
\prod_{i\in S} \tau_{k_i}(1) \right|
\nu\rang\,,
\end{multline}
where the summation is over all subsets $S$ and all intermediate
partitions $\eta$ of the same size as $\mu$ and $\nu$. 

After repeated application, the recursive
relation \eqref{rcr} expresses the original rubber integrals
\eqref{ep1} in terms of two types of integrals: 
rubber integrals {\em without} markings and tube integrals
without cotangent line classes at the relative points.

The following Lemma completes the 
reduction of the rubber integrals \eqref{rcr} to
the tube theory \eqref{ttt}.

\begin{lm}
\label{mt}
For $n=0$,
\begin{equation}
  \label{rbh}
  \lang \mu,k \left| \right|
\nu\rang^\sim = \frac{1}{(k+1)!} 
\lang \mu \left| \tau_1(\omega)^{k+1} \right| \nu \rang  \,.
\end{equation}
\end{lm}

\bpf
The dilaton equation for the relative theory yields: 
\begin{eqnarray*}
\lang \mu,k \left|  \psi_1 \right| \nu 
\rang_g^\sim & = & (2g-2+\ell(\mu)+\ell(\nu)) 
\lang \mu,k \left| \right| \nu 
\rang_g^\sim \\
& = & (k+1) \lang \mu,k \left| \right| \nu \rang_g^{\sim}
\end{eqnarray*}
The second equality is obtained by matching $k$ with the
dimension of the moduli space $R^\bullet_{g,0}(\mu,\nu)$.
If $k=0$, the Lemma is proven by the dilaton equation and
Lemma \ref{sittt}.

We proceed by induction.
For $k>0$, the induction step is proven
by
applying the recursion \eqref{rcr} to the left side
of the dilaton equation
and using the degeneration formula.
\epf

The reduction of the rubber integrals to the tube theory will be
expressed in terms of the operator formalism of the
infinite wedge representation in Section \ref{two}.
The following convention will simplify our
formulas:
\begin{equation}
  \label{rbh-1}
  \lang \mu,-1 \left| \right|
\nu\rang^\sim = 
\lang \mu \left|\right| \nu \rang  = 
\frac{\delta_{\mu,\nu}}{\zz(\mu)}\,.
\end{equation}
For example, the recursion \eqref{rcr} is valid for $k=0$ with
the above convention.

\subsubsection{}
The outcome of the rubber calculus is a determination of the
depth $n$ rubber integral \eqref{ep1} in terms of the depth $n-1$
theory of the tube.

\begin{pr}\label{tcap}
The depth $r$  theory of the tube determines the depth $r+1$
theory of the cap.
\end{pr}

\bpf
By Proposition \ref{lcapp},
the equivariant theory of the cap is determined by
the localization formula.
After nonequivariant specialization, the depth $r+1$ theory
of the cap is determined by depth $r+1$ rubber integrals.
By the rubber calculus, the latter integrals are determined by
the depth $r$  theory of the tube.
\epf

\subsection{Virasoro for the tube}

The Virasoro constraints for the tube take a very simple form.
Since $$\chi({\mathbf T}^*)=0,$$ there are no extensive Virasoro
reactions.

\begin{pr}
\label{virtube}
  The Virasoro constraints for the theory of
the tube are equivalent to the following 
relation: 
  \begin{equation}
    \label{tvi}
    \lang \mu \left| \tau_{l}(\omega) \, 
\prod \tau_{k_i}(1) \right|
\nu\rang = \binom{\sum k_i + l}{k_1,\dots,k_n,l}  
\,
\lang
\mu \Big| \tau_{l+\sum (k_i-1)}(\omega)\Big| \nu \rang \,.
  \end{equation}
\end{pr}

\begin{proof}
Equality \eqref{tvi} 
follows from the Virasoro constraints 
by induction using the following elementary identity for
multinomial coefficients:
\begin{multline*}
  \binom{\sum k_i + l}{k_1,\dots,k_n,l} = 
\binom{k_1 + l}{k_1} \, \binom{\sum k_i + l-1}{k_2,\dots,k_n,k_1 + l-1} + \\
\sum_{j=2}^n \binom{k_1 + k_j-1}{k_1} \,
\binom{\sum k_i + l-1}{k_2,\dots,k_j+k_1-1,\dots,k_n,l} \,.
\end{multline*}
Indeed, relation \eqref{tvi} is equivalent to validity of the
Virasoro constraints for relative invariants of the tube
with exactly 1 stationary descendent.

Relation \eqref{tvi}, together with the 
degeneration formula, uniquely determines the entire relative theory,
\begin{equation}
\label{fulll}
\lang \mu \left| \prod_{i} \tau_{k_i}(\omega) \prod_j \tau_{l_j}(1)
\right| \nu \rang,
\end{equation}
in terms of the stationary theory:
\begin{enumerate}
\item[(i)]
If there are no stationary descendents, then integral \eqref{fulll}
vanishes. Proofs of the vanishing are easily obtained by degeneration
arguments or the
localization formula for the tube.  
\item[(ii)]
If there is at least one stationary descendent,
we can degenerate the tube $\mathbf T$ into a chain of tubes,
each containing exactly one stationary descendent.
\end{enumerate}
Since the Virasoro constraints are compatible with degeneration,
relation \eqref{tvi} implies the validity of the Virasoro constraints
for the entire theory of the tube.
\end{proof}
\vspace{+10pt}

We will only require the nonequivariant generating series for
the tube:
\begin{multline*}
  \bG'(\mu \left|x_1,\dots,x_n,y_1,\dots,y_m\right|\nu) = \\
\sum_g \sum_{k_i\ge 0} \sum_{l_j\ge 0} \, \prod_{i=1}^n x_i^{k_i+1} 
\prod_{j=1}^m y_j^{l_j+1} \, 
\lang \mu \left| \prod \tau_{k_i}(\omega)  \, \prod \tau_{l_j}(1) 
\right| \nu \rang_g  \,.
\end{multline*}
Relation 
\eqref{tvi} is equivalent to the following equation:
\begin{multline}
  \label{tvi2}
  \bG'(\mu|x,y_1,\dots,y_m|\nu) = \\
xy_1\cdots y_m \left(x+\sum y_j\right)^{m-1} \, 
\bG'\left(\mu\left|x+\sum y_j\right|\nu\right) \,.
\end{multline}
The function 
\begin{align}
  \bT(x_1,\dots,x_k) &= x_1 \dots x_k \, \left(\sum x_i\right)^{k-2}
\notag \\
& = \sum_{T} \prod x_i^{\val_T(i)} \label{sumT} 
\end{align}
will occur often. The sum in \eqref{sumT} is over all trees $T$ with 
vertex set $\{1,\dots,k\}$. The valence of the vertex $i$ of $T$
is $\val_T(i)$.
The factorization of the tree sum is
a classical result due to Cayley.

\subsection{Outlook}
By Proposition \ref{fff}, the depth $r$ theory of the cap determines the 
depth $r$ theory of the tube. 
By Proposition \ref{tcap}, the depth $r$  theory of the tube
determines the depth $r+1$  theory of the cap.
Together, the results uniquely
determine the entire relative theories of both the cap and the tube from
their stationary theories.
We prove the Virasoro constraints for ${\mathbf C}$ and 
$\mathbf T$  by showing the
Virasoro operators are compatible with the
opposite determinations.

The proof of the Virasoro constraints for the full 
relative theory
of target curves requires several additional techniques for handling the
odd cohomology.

\pagebreak

\section{The operator formalism}
\label{oopp}

\subsection{Review of the infinite wedge space}

\label{iwr}

\subsubsection{}

Let $V$ be a linear space with basis $\left\{\ul{k}\right\}$
indexed by the half-integers:
$$V = \bigoplus_{k \in \Z+ \sh} \com \, \ul{k}.$$
For each subset $S=\{s_1>s_2>s_3>\dots\}\subset \Z+\sh$ satisfying:
\begin{enumerate}
\item[(i)] $S_+ = S \setminus \left(\Z_{\le 0} - \sh\right)$ is
finite,
\item[(ii)]$S_- =
\left(\Z_{\le 0} - \sh\right) \setminus S$ is finite,
\end{enumerate}
we denote by  $v_S$ the following infinite wedge product:
\begin{equation}
v_S=\ul{s_1} \wedge \ul{s_2} \wedge  \ul{s_3} \wedge \dots\, .
\label{vS}
\end{equation}
By definition,
 $$\LV= \bigoplus \com \,v_S
$$
is the linear space with basis $\{v_S\}$.
Let $(\,\cdot\,,\,\cdot\, )$ be the inner product on  $\LV$ for which
$\{v_S\}$ is an orthonormal basis.

\subsubsection{}

The fermionic operator $\psi_k$ on $\LV$ is defined by wedge product with
the vector $\ul{k}$,
$$
\psi_k \cdot v = \ul{k} \wedge v  \,.
$$
The operator $\psi_k^*$ is defined as the adjoint of $\psi_k$
with respect to the inner product  $(\,\cdot\,,\,\cdot\, )$.

These operators satisfy the canonical anti-commutation relations:
\begin{gather}
  \psi_i \psi^*_j + \psi^*_i \psi_j = \delta_{ij}\,, \\
  \psi_i \psi_j + \psi_j \psi_1 = \psi_i^* \psi_j^* + \psi_j^*
  \psi_i^*=0.
\end{gather}
The {\em normally
ordered} products are defined by:
\begin{equation}\label{e112}
\nr{\psi_i\, \psi^*_j} =
\begin{cases}
\psi_i\, \psi^*_j\,, & j>0 \,,\\
-\psi^*_j\, \psi_i\,, & j<0 \,.
\end{cases}
\end{equation}

\subsubsection{} \label{sCH}

Let $E_{ij}$, for $i,j\in \Z+\sh$, be the standard basis of
matrix units of $\gli$.
The assignment
$$
E_{ij} \mapsto\, \nr{\psi_i\, \psi^*_j}\ \ ,
$$
defines
a projective representation of the Lie algebra $\gli=\mathfrak{gl}(V)$ on
$\LV$.

The {\em charge} operator $C$ corresponding to the identity matrix of $\gli$,
$$
C= \sum_{k\in\Z+\frac12} \, E_{kk},
$$
acts on the basis $v_S$ by:
$$
C \, v_S = (|S_+| - |S_-|) v_S \,.
$$
The kernel of $C$, the zero charge subspace,
is spanned by the vectors
\begin{equation*}
v_{\lam}=\ul{\lam_1-\tfrac12} \wedge \ul{\lam_2-\tfrac32} \wedge
\ul{\lam_3-\tfrac52} \wedge \dots
\end{equation*}
indexed by all partitions $\la$. We will denote the kernel by $\LVc$.

The eigenvalues on $\LVc$ of the {\em energy} operator,
$$
H= \sum_{k\in\Z+\frac12} k \, E_{kk},
$$
are easily identified:
$$
H \, v_\lambda = |\lambda| \, v_\lambda\,.
$$
The vacuum vector
$$
\vac = \ul{-\tfrac12} \wedge \ul{-\tfrac32} \wedge \ul{-\tfrac52}
\wedge \dots
$$
is the unique vector with the minimal (zero) eigenvalue of $H$.

The {\em vacuum expectation} $\lang A \rang$ of an operator $A$ on $\LV$ is defined
by  the inner product:
$$\lang A \rang = (A v_\emptyset, v_\emptyset).$$

\subsubsection{}

For any $r\in\Z$, we define
\begin{equation}
\cE_r(z) = \sum_{k\in\Z+\frac12} \, e^{z(k-\frac{r}2)} \,  E_{k-r,k}
+ \frac{\delta_{r,0}}{\cs(z)} \, ,\label{Er}
\end{equation}
where the function $\cs(z)$ is defined by
\begin{equation}
  \label{cs}
  \cs(z) = e^{z/2} - e^{-z/2} \,.
\end{equation}
The exponent in \eqref{Er} is set to satisfy:
$$
\cE_r(z)^* = \cE_{-r}(z)\,,
$$
where the adjoint is with respect to the standard inner
product on $\LV$.

Define the operators $\cP_k$ for $k>0$
by:
\begin{equation}
\frac{\cP_k}{k!} = \, [z^k] \, \cE_0(z)\,,\label{defcP}
\end{equation}
where $[z^k]$ stands for the coefficient of $z^k$.
The operator,
$${\mathcal F}_2 = \frac{\cP_2}{2!}= \sum_{k\in\Z+\frac12} \frac{k^2}{2} E_{k,k}\, ,$$
will play a special role.

\subsubsection{}

The operators $\cE$ satisfy the following fundamental
commutation relation:

\begin{equation}\label{commcEE}
  \left[\cE_a(z),\cE_b(w)\right] =
\cs\left(\det  \left[
\begin{smallmatrix}
  a & z \\
b & w
\end{smallmatrix}\right]\right)
\,
\cE_{a+b}(z+w)\,.
\end{equation}

Equation \eqref{commcEE}
automatically incorporates the central
extension of the $\gli$-action, which appears as the
constant term in $\cE_0$ when $r=-s$.

\subsubsection{}

The operators $\cE$ specialize to the
standard bosonic operators on $\LV$:
$$
\al_k = \cE_k(0)\,, \quad k\ne 0 \,.
$$
The commutation relation \eqref{commcEE} specializes to
the following equation
\begin{equation}
\label{qhhqq}
[\al_k, \cE_l(z)] =
 \cs(kz) \, \cE_{k+l}(z) \,.
\end{equation}
When $k+l=0$, equation \eqref{qhhqq} has the following constant term:
$$
\frac{\cs(kz)}{\cs(z)}=\frac{e^{kz/2}-e^{-kz/2}}{e^{z/2}-e^{-z/2}} \,.
$$
Letting $z\to 0$, we recover the standard relation:
\begin{equation*}
[\al_k,\al_l]= k \, \delta_{k+l} \,.
\end{equation*}

\subsection{Rubber calculus II}

\subsubsection{}
We will find an operator formula for rubber integrals
using the calculus of Section \ref{pu} and
Proposition \ref{virtube}.
In particular, the operator formula here {\em will
depend upon the Virasoro constraints for the tube}.

\subsubsection{}
For a partition $\nu$, let the vector 
$\left.|\nu\rang\in \LV$ be defined by 
$$
\left|\nu\rang
 = \frac1{\zz(\nu)} \, \prod \al_{-\nu_i} \, \vac \,.
$$
The operator
$$
\sum_{|\eta|=d} \Big| \eta \Big\rangle \, \zz(\eta)\, \Big\langle
 \eta \Big|
$$
is the orthogonal projection onto the subspace of the 
Fock space $\LV_0$ corresponding to partitions of size $d$.

\subsubsection{}

The $1$-point stationary generating series of the tube
has been calculated in \cite{OP2} by the GW/H correspondence:
$$
\bG'\left(\mu\left|x\right|\nu\right) = 
\lang \mu \left|\cE_0(x) \right|
\nu \rang \,.
$$
Introduce the following operator,
$$
\cE_0(x_1,\dots,x_k) = \bT(x_1,\dots,x_k) \, \cE_0(x_1+\dots+x_k) \,.
$$
Then, equation \eqref{tvi2} yields:
\begin{equation}
  \label{tvi3}
  \bG'(\mu|x,y_1,\dots,y_m|\nu) = 
\lang \mu \left|\cE_0(x,y_1,\dots,y_m) \right|
\nu \rang  \,. 
\end{equation}

\subsubsection{}

Similarly, by  Lemma \eqref{mt}, 
\begin{equation}
\sum_{k\ge -1} s^{k+1} \, \lang \mu,k \left| \right|
\nu\rang^\sim = \sum_{k\ge -1} \frac{s^{k+1}}{(k+1)!} \, 
\lang \mu \left| {\mathcal F}_2^{k+1} \right| \nu \rang = 
\lang \mu \left| e^{s {\mathcal F}_2} \right| \nu \rang \,.
\label{tvi4}
\end{equation}
Here, we follow convention \eqref{rbh-1}.

\subsubsection{}


Let $\bE(x_1,\dots,x_m,s)$
denote the following operator:
\begin{equation}
  \label{bE}
  \bE(x_1,\dots,x_m,s) = \sum_{\pi} s^{m-\ell(\pi)} 
\prod_{k=1}^{\ell(\pi)} \cE_0(x_{\pi_k}) \,,
\end{equation}
where the summation is over all partitions 
$$
\pi= \pi_1 \sqcup \dots \sqcup \pi_{\ell}
$$
of the set $\{1,\dots,m\}$ into nonempty disjoint
subsets. Here, $x_{\al_k}$ denotes the variables $x_i$
with indices in the subset $\pi_k$. For example
\begin{multline*}
  \bE(x_1,x_2,x_3,s)= \cE_0(x_1) \, \cE_0(x_2) \, \cE_0(x_3) + 
s \, \cE_0(x_1,x_2) \, \cE_0(x_3)\\
s \, \cE_0(x_1,x_3) \, \cE_0(x_2) +
s \, \cE_0(x_2,x_3) \, \cE_0(x_1)  + s^2 \, \cE_0(x_1,x_2,x_3) \,.
\end{multline*}
Since the operators $\cE_0$ commute, the ordering of
the blocks of $\pi$ is not important. 

The recursive relation
\eqref{rcr}  together with formulas \eqref{tvi3} and \eqref{tvi4}
directly yields an operator formula for the generating
series $\bG^\sim$ of rubber integrals.

\begin{pr}
\label{Gto}
\begin{equation*}
    \bG^\sim(\mu,s|x_1,\dots,x_m|\nu) = 
\lang \mu \left| e^{s \cP_2} \bE(x_1,\dots,x_m,\tfrac 1s) \right|
\nu \rang  \,.
\end{equation*}
\end{pr}

The validity of Proposition \ref{Gto} depends upon the Virasoro constraints
for the {depth} $m-1$ theory of the tube.

\subsection{The equivariant relative cap revisited}

Recall the basic operators $\bA(z)$ defined in \cite{OP3}:
\begin{equation}\label{bA}
   \bA(z) = 
\cS(z)^{t z} \, 
\left(\sum_{k\in \Z}
\frac{\cs(z)^k}{(tz+1)_k } \, \cE_k(z) 
\right) \, ,
\end{equation}
where
$$
\cs(z)=e^{z/2}-e^{-z/2}\,, \quad \cS(z)=\frac{\cs(z)}{z} =
\frac{\sinh{z/2}}{z/2} \,
$$
and
$$
(a+1)_k = \frac{(a+k)!}{a!} =
\begin{cases}
 (a+1) (a+2) \cdots (a+k) \,, & k\ge 0 \,,  \\
 (a (a-1) \cdots (a+k+1))^{-1}  \,, & k \le 0 \,.
\end{cases}
$$
Here, the genus variable $u$ of \cite{OP3} is set to 1.

Proposition \ref{lcapp} together with the operator formulas for $\bH$ of
 Section 5 of \cite{OP3} and the operator formula for $\bG^\sim$
yields:
\begin{multline}
  \label{loc1}
  \bG(z_1,\dots,z_n,w_1,\dots,w_m|\nu) = \\
\sum_\eta \lang \left.
\prod \bA(z_i) \, e^{\alpha_1} \, e^{\frac1t \, {\mathcal F}_2} 
\right| \eta \rang \, \zz(\eta)\, \lang \eta \left|
e^{-\frac1t {\mathcal F}_2} \, \bE(w_1,\dots,w_m,-t) \right|
\nu \rang \,,
\end{multline}
where the summation is over all partitions $\eta$ such 
that $|\eta|=|\nu|$.

The orthogonal projection,
$$
\sum_{|\eta|=d} \Big| \eta \Big\rangle \, \zz(\eta)\, \Big\langle
 \eta \Big|
$$
commutes with ${\mathcal F}_2$ and $\bE$ and fixes the vector 
$\left|\nu\rang$. After commuting the projection  to the far
right, we obtain the following the operator formula for the
equivariant relative theory of the cap.

\begin{pr} \label{lcapp2}
\begin{multline}
  \label{loc}
  \bG(z_1,\dots,z_n,w_1,\dots,w_m|\nu) = \\
\lang \left.
\prod \bA(z_i) \, e^{\alpha_1} \, 
 \bE(w_1,\dots,w_m,-t) \right|
\nu \rang \,.
\end{multline}
\end{pr}

The validity of Proposition \ref{lcapp2} depends upon the Virasoro constraints
for the {depth} $m-1$ theory of the tube.

\pagebreak

\section{Virasoro constraints for the cap}
\label{llee}

\subsection{Plan of the proof}
We prove the Virasoro constraints for the cap by induction on depth.
The base of the induction is trivial as 
the Virasoro constraints are empty for the depth 0 theory of the cap.
Assume the Virasoro constraints hold for the 
depth $r$ theory of the cap. By Proposition \ref{eee}, the Virasoro constraints
hold for the depth $r$ theory of the tube. Then, 
Proposition \ref{lcapp2} provides operator formulas for the equivariant
series,
\begin{equation}
\label{es}
\bG(z_1,\dots,z_n,w_1,\dots,w_m|\nu),
\end{equation}
for $m\leq r+1$.
Via the nonequivariant limit, the series \eqref{es} determine
the depth $r+1$ nonequivariant theory of the disk.
To complete the induction step, we must prove the
nonequivariant limits of the operator formulas for
\eqref{es} satisfy the Virasoro rules.

\subsection{The nonequivariant limit}

\subsubsection{} 
Our goal now is to study the nonequivariant limit
of the operator formula \eqref{loc} for 
$$\bG(z_1,\dots,z_n,w_1,\dots,w_m|\nu).$$
We will obtain an 
operator formula for the nonequivariant limit  involving
the first two terms,  $\bA^0(z)$ and $\bA^1(z)$, in the expansion 
of the operator $\bA(z)$,
$$
\bA(z) = \sum_{k\ge 0} t^k \, \bA^k(z), 
$$
in powers of the parameter $t$. 
The operator 
$$
\bA^0(z) = \bA(z)\big|_{t=0}
$$
is related to the operator $\cE_0(z)$ via
\begin{equation}
e^{-\alpha_1}\, \bA^0(z) \, e^{\alpha_1} = \cE_0(z)  \, ,
\label{AaE}
\end{equation}
see \cite{OP3}.
The above relation will play an important role in our
analysis of formula \eqref{loc}\,.

\subsubsection{}
\label{veee}

We will often encounter the operator $\bA^\vee$ constructed 
from
the operators $\bA^0$ and $\bA^1$:
\begin{multline}\label{Avee}
 \bA^\vee(z_1,\dots,z_n) = \bT(z_1,\dots,z_n) \, \Bigg[ \\
z_n \, \ln\left(1+\frac{z_1+\dots+z_{n-1}}{z_n}\right) \,
\bA^0(z_1+\dots+z_n)  + \\
\frac{z_n}{z_1+\dots+z_n} \, \bA^1(z_1+\dots+z_n) \Bigg] \,.
\end{multline}
The expression in square brackets is the coefficient
of $t$ in the expansion of 
$$
\left(\frac{z_1+\dots+z_{n}}{z_n}\right)^{tz_n} \,
\bA(tz_n,z_1+\dots+z_n) 
$$
in powers of $t$, where the $2$-parameter operator
$\bA(a,b)$ is defined by
\begin{equation}\label{bAab}
   \bA(a,b) = 
\cS(b)^{a} \, 
\left(\sum_{k\in \Z}
\frac{\cs(b)^k}{(a+1)_k} \, \cE_k(b) 
\right) \,.
\end{equation}

\subsubsection{}

An operator formula for 
the nonequivariant limit of $\bG(z,w)$ is
obtained from Proposition \ref{lcapp2} and the relation,
$$1=\frac{[0]-[\infty]}{t},$$
in the localized equivariant cohomology of $\proj^1$,
\begin{equation}
  \label{nloc}
  \bG'(z_1,\dots,z_n,w_1,\dots,w_m|\nu) = [t^m] \, 
\lang \left.
\prod \bA(z_i) \, \bM(w_1,\dots,w_m) \right|
\nu \rang\,,
\end{equation}
where $[t^m]$ denotes the coefficient of $t^m$.  The operator 
$\bM(w)$ is defined by:
\begin{equation}
  \label{bM}
  \bM(w_1,\dots,w_m) = \sum_{S\subset \{1,\dots, m\}}
(-1)^{|S|} \prod_{i\notin S} \bA(w_i) \, e^{\alpha_1} \, 
 \bE(w_S,-t) \,.
\end{equation}
where the summation is over all subsets $S$ of the 
index set $\{1,\dots, m\}$. Here, $w_S$ denotes the 
variables $w_j$ such that $j\in S$. 

Inside the brackets, the order of the operators in the product
$$\prod_{i\notin S} \bA(w_i)$$
does not matter by the symmetry of the
equivariant series $\bG$. However, the
operators themselves do {\em not} commute.
We order the operators in the above product by their indices $i$.

\subsubsection{}
We now express the leading $t$ coefficient of the
operator $\bM$ in terms of the operators
$\bA^\vee$ introduced in Section \ref{veee}.

\begin{pr}
\label{tM}
We have

\begin{enumerate}
\item[(i)]
for $k <m$,
\begin{equation}
  \label{vanc}
   [t^k] \, \bM(w_1,\dots,w_m) =0\,, 
\end{equation}
\item[(ii)] for $k=m$,
\begin{equation} \label{Ml}
    [t^m] \, \bM(w_1,\dots,w_m) = \sum_{\pi} 
\left(\timeprod_{p=1,\dots,\ell(\pi)} \bA^\vee (w_{\pi_p}) 
\right)\, e^{\alpha_1} \,.
\end{equation} 
\end{enumerate}
\end{pr}

The summation in \eqref{Ml} is over all partitions $\pi$ of the 
 ordered set $\{1,\dots, m\}$.
The argument $w_{\pi_p}$ denotes the set of 
variables in the induced order with indices in the part $\pi_p$ of $\pi$. 
The
operators $\bA^\vee$ are ordered left-to-right by increasing
maximal elements of $\pi_p$. 
Since the operators $\bA^\vee$ do not commute,  their ordering is
important.

\subsubsection{}
The operator $\bA(z)$ is obtained from $\bA(a,b)$ by the
evaluation
$$
\bA(z) = \bA(tz,z) \,.
$$
Therefore,
$$
\bA^0(z) = \bA(0,z)\,.
$$
The following result will play an important role in the proof
of Proposition \ref{tM}. 

\begin{lm} We have \label{bbAAcc}
\begin{equation}\label{bAc}
\left[\bA^0(z), \bA(a,b)\right] = 
\cs(za)\, \left(1+\frac{z}{b}\right)^a \, \bA(a,z+b) \,.
\end{equation}
\end{lm}

\begin{proof}
The basic commutation relation for the operators $\cE$ is:
\begin{equation}
  \label{cmE}
 \left[\cE_k(z),\cE_l(b)\right] = \cs(kb-lz) \, \cE_{k+l}(z+b) \,. 
\end{equation}
Therefore,
\begin{equation*}
  \left[\bA^0(z), \bA(a,b)\right] = 
\cS(b)^a \,\sum_{k\ge 0,\,l\in \Z} \, 
\frac{\cs(z)^k \, \cs(b)^l \, \cs(kb-lz)}
{k! \, (a+1)_l } \,
\cE_{k+l}(z+b) \,.
\end{equation*}
The coefficient of $\cE_{m}(z+b)$ in the above series equals
\begin{equation}
\cS(b)^a \, \frac{\cs(b)^m}{(a+1)_m}\, 
\sum_{k\ge 0} \binom{m+a}{k} \, \cs(kb+kz-m z)\, 
\left(\frac{\cs(z)}{\cs(b)}\right)^k \,.
\label{cfEm}
\end{equation}
Since
$$
\cs(kb+kz-m z)\, 
\left(\frac{\cs(z)}{\cs(b)}\right)^k =
e^{-mz/2} \, \left(\frac{e^z-1}{1-e^{-b}}\right)^k  -
e^{mz/2} \, \left(\frac{1-e^{-z}}{e^{b}-1}\right)^k \,,
$$
the binomial series in \eqref{cfEm} sums to 
$$
e^{-mz/2}\,\left(\frac{e^z-e^{-b}}{1-e^{-b}}\right)^{m+a} -
e^{mz/2}\,\left(\frac{e^{b}-e^{-z}}{e^{b}-1}\right)^{m+a} =
\cs(za) \, \frac{\cs(b+z)^{m+a}}{\cs(b)^{m+a}} \,.
$$
Thus, the coefficient  \eqref{cfEm} of $\cE_m(z+b)$ in the 
commutator equals
$$
\frac{\cs(za)}{b^a} \, \frac{\cs(z+b)^{m+a}}{(a+1)_m} \,,
$$
as was to be shown. 
\end{proof}


\begin{cor} We have
\begin{multline}\label{bAcm}
\left[\bA^0(z_1),\left[\bA^0(z_2),
\left[\dots,\left[\bA^0(z_n),\bA(a,b)\right]\dots\right]\right]\right]
 = \\
\prod \cs(z_ia)\, \left(1+\frac{\sum z_i}{b}\right)^a \, 
\bA\left(a,\sum z_i+b
\right) \,.
\end{multline}
\end{cor}

\subsubsection{}\label{alg}

The strategy of our proof 
of Proposition \ref{tM} will be presented  here.
The proof is given in Sections \ref{azaa} -- \ref{bzbb}.

We must extract  coefficients
of $t$ from the operator $\bM$ defined by formula \eqref{bM}.
Consider the summand of the formula indexed by $S$:
\begin{equation}
\label{jj}
(-1)^{|S|} \prod_{i\notin S} \bA(w_i) \, e^{\alpha_1} \, 
 \bE(w_S,-t).
\end{equation}
Let $\bA(w_i)$ be the leftmost operator in \eqref{jj}
where $$
i=\min\{j,j\in \CS\}\,,
$$
and $\CS$
denotes the complement of $S$. 
We start by
selecting a coefficient of $t$ from $\bA(w_i)$. 
There are two possibilities. 
Either we select the minimal power $t^0$, 
corresponding to the operator  $\bA^0(w_i)$, or we select {\em all}
higher power of $t$: 

\begin{enumerate}
\item[(i)] If we select the minimal
power $t^0$ in $\bA(w_i)$, we commute the resulting operator  $\bA^0(w_i)$
to the right
of the operator $e^{\alpha_1}$.
The operator $\bA^0(w_i)$ commutes through the 
operators $\bA(w_j)$ via the commutation relation 
\begin{equation}\label{bAw}
\left[\bA^0(w_i), \bA(w_j)\right] = 
\cs(t w_i w_j)\, \left(1+\frac{w_i}{w_j}\right)^{tw_j} \, 
\bA(tw_j,w_i+w_j) \,,
\end{equation}
obtained from Lemma \ref{bbAAcc}. 

The
prefactors in equation \eqref{bAw} play an
important role.
In particular, since
$$
\cs(t w_i w_j) = t w_i w_j + o(t^2) ,
$$
the minimal positive power of $t$ that can be extracted from the 
operator \eqref{bAw} is $t^1$, which yields the operator
$w_i w_j \, \bA^0(w_i+w_j)$. 

When the commutator \eqref{bAw}
appears, we say the $i$th operator interacts with 
the $j$th operator. 
The final commutation of $\bA(w_i)$ with $e^{\alpha_1}$ is determined
by relation \eqref{AaE}. When the relation \eqref{AaE} is
used, we say the $i$th operator interacts with $\infty$.

\item[(ii)]
If the higher powers of $t$ in $\bA(w_i)$ are chosen at the start, then we do
nothing further with $\bA(w_i)$.
\end{enumerate}

Once step (i) or (ii) has been completed, we proceed to treat the 
leftmost operators of the remaining terms, repeating the 
above cycle  using step (i) for the minimal power of $t$ and
step (ii) for the higher powers. 
The cycle is repeated until all the operators $\bA$
of all the terms are treated.
The outcome is a transformation of formula
\eqref{bM}.
Proposition \ref{tM} is proven by finding 
 a large cancellation in the
transformed formula.

\subsubsection{} \label{azaa}

\begin{figure}[htbp]\psset{unit=0.5cm}
  \begin{center}
    \begin{pspicture}(0,-1)(20,5)
\psarc{<-}(4,0.2){4}{0}{180}
\psarc{<-}(7,0.2){1}{0}{180}
\psarc{<-}(7,0.2){3}{0}{180}
\psarc{<-}(12,0.2){2}{0}{180}
\psarc{<-}(16,0.2){4}{0}{180}
\psdots[dotscale=1](0,0)(2,0)(4,0)(6,0)(8,0)(10,0)(12,0)(14,0)
\rput[l](0.4,0){$i$}
\rput[l](8.4,0){$j$}
\rput[t](20,0){$\infty$}
\rput[t](7,-0.5){$\underset{\CS}{\underbrace{\hspace{7.3cm}}}$}
\end{pspicture}
   \caption{An interaction diagram}
    \label{fint}
  \end{center}
\end{figure}
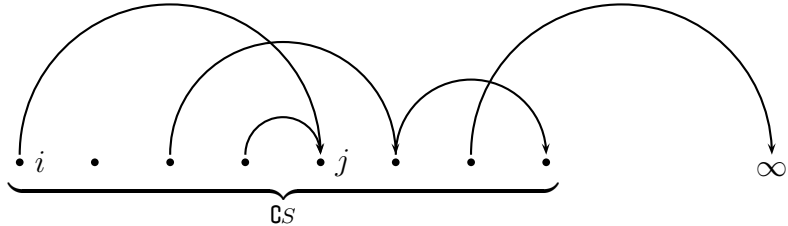

The following diagrammatic technique is useful for
understanding the terms arising in this 
transformation of the formula \eqref{jj}. The terms
of the transformed formula will be indexed by
graphs $\cI$ of the kind shown in Fig.~\ref{fint} which 
we call \emph{interaction diagrams}. 
An interaction diagram has vertices indexed by the set $\CS$ plus 
an additional vertex marked by the symbol $\infty$. 
If the $i$th operator interacts with the $j$th 
operator, we draw an arrow from $i$ to $j$. The notation 
$i\to j$ will also be used. 
If the $i$th operator interacts with $\infty$, we draw
an arrow from $i$ to $\infty$. A vertex with no
outgoing arrows 
in an interaction diagram marks an occurrence of
the step (ii) above. 

Abstractly, interaction diagrams $\cI$ can be defined
as trees with vertex set $\CS\cup\{\infty\}$ such that
for every $i\in \CS$ there is at most one edge 
$(i,j)$ with $j>i$. 
An interaction diagram $\cI$ gives the set $\CS$
a natural partial order $\ioe$ in which by
definition $i\ioe j$ if
$i$ can be connected to $j$ by arrows of $\cI$. 
Note that this order is compatible with the
usual ordering of natural numbers.

\subsubsection{}

For any $k\in \CS$, let $\cI(k)$ be the maximal 
connected subdiagram of $\cI$ with maximal vertex $k$. 
Clearly, 
\begin{equation}
  \label{vertJ}
  \vrt(\cI(k)) = \{i \, | \, i \ioe k\} \,. 
\end{equation}
The operators $\bA(w_i)$, $i\ioe k$, are
precisely the operators that influence the 
operator $\bA(w_k)$ in the transformation of
the term \eqref{jj}. 

Concretely, the operator $\bA(w_k)$ is transformed 
by the interactions specified
by $\cI(k)$ to
the 
operator: 
\begin{equation}
  \label{bbA}
  \ww(\cI(k)) \, \left(1+\frac{\sum_{i\io k} w_{i}}{w_k}\right)^{tw_k} \, 
\bA\left(t w_k,\sum_{i\ioe k} w_{i}
\right)\,,
\end{equation}
where $\ww(\cI(k))$ is the following product over all 
arrows $i\to j$ in  $\cI(k)$ 
$$
\ww(\cI(k)) = \prod_{i\to j} 
\cs\left(t\left(\sum_{i'\ioe i} w_{i'}
\right)w_j\right)\,.
$$
In \eqref{bbA}, only the factor $\ww(\cI(k))$  depends on 
the actual interaction diagram $\cI(k)$ and not only on the 
vertex set \eqref{vertJ}. Therefore, it is natural to sum the weights
$\ww(\cI(k))$ over all interaction diagram with the same
vertex set $K$, where 
$$
K=\{ k_1 < k_2 < \dots < k_r=k\} \,,
$$
is an arbitrary subset of $\CS$ with maximal element $k$.

\begin{lm}\label{sspp}
We have, 
\begin{equation}
  \label{smI}
  \sum_{\vrt(\cI(k))=K} \, \ww(\cI(k)) =
t^{r-1} \, \bT(w_{k_1},\dots,w_{k_r}) + o(t^r)\,,
\end{equation}
where $\bT$ is the tree function defined by \eqref{sumT}.
\end{lm}

\bpf
We present a combinatorial proof of the following
equivalent identity:
\begin{equation}
  \label{oT}
  \sum_{\vrt(\cI(k))=K} \,\,
\prod_{i\to j} \left(\sum_{i'\ioe i} w_{i'}
\right)w_j = 
\bT(w_{k_1},\dots,w_{k_r}) \,. 
\end{equation}
We will construct a
map $\bo$ from the set of trees $T$ on $K$
to the set of interaction diagrams $\cI(k)$ with the
following property:
\begin{equation}
\prod_{i\to j} \left(\sum_{i'\ioe i} w_{i'}
\right)w_j =
\sum_{T\in \bo^{-1}(\cI(k))} \prod w_i^{\val_T(i)}
\,.
\label{TO}
\end{equation}
Identity \eqref{oT} is then deduced from the definition of $\bT$.

Let $T$ be an arbitrary tree on the vertex set $K$ with
edge set $E$. 
What needs to be constructed is the edge set of 
the interaction diagram $\sigma(T)$. 
Denote this edge set by $\sigma(E)$. The
construction of $\sigma(E)$ will be inductive. 

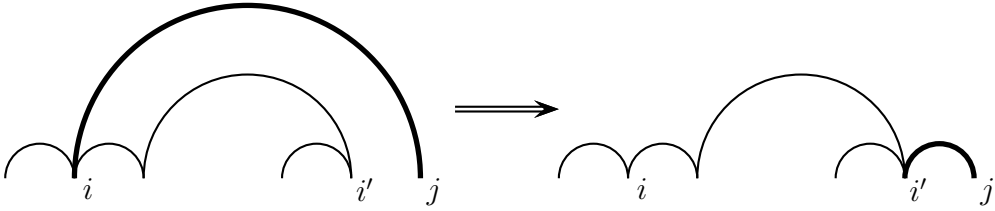
\begin{figure}[htbp]\psset{unit=0.46cm}
  \begin{center}
    \begin{pspicture}(0,-1)(30,5)
\psline[doubleline=true]{->}(13,2)(16,2)

\psarc(1,0){1}{0}{180}
\psarc(3,0){1}{0}{180}
\psarc(7,0){3}{0}{180}
\psarc(9,0){1}{0}{180}
\psarc[linewidth=2pt](7,0){5}{0}{180}

\rput[t](2.4,0){$i$}
\rput[t](10.4,0){$i'$}
\rput[t](12.4,0){$j$}

\psarc(17,0){1}{0}{180}
\psarc(19,0){1}{0}{180}
\psarc(23,0){3}{0}{180}
\psarc(25,0){1}{0}{180}
\psarc[linewidth=2pt](27,0){1}{0}{180}

\rput[t](18.4,0){$i$}
\rput[t](26.4,0){$i'$}
\rput[t](28.4,0){$j$}

\end{pspicture}
   \caption{From trees to interaction diagrams}
    \label{fcs}
  \end{center}
\end{figure}

 Let $(i,j)\in E$
be an edge of $T$ and assume that all edges $(a,b)\in E$
with $b<j$ have been already considered. Then we add 
to $\sigma(E)$ the edge $i'\to j$, where $i'$ is the maximal vertex connected to $i$
by edges already in $\sigma(E)$. In other words, 
the new edge is obtained from the edge $(i,j)$ by 
moving the origin as far as we can along the edges
already in $\sigma(E)$. This is illustrated in Fig.~\ref{fcs}. 

The (multivalued) inverse of the map $\bo$ is 
constructed similarly.  The edge set $E$ of a tree 
 $T\in \bo^{-1}(\cI(k))$ is be constructed inductively. Suppose all 
edges with endpoints less than $j\in K$
have already been constructed. Let $i'\to j$ be an
arrow in $\cI(k)$ with endpoint $j$. Then, we add to $E$
an arbitrary edge $(i,j)$ where the vertex $i$
is connected to $i'$ by the edges already in $E$.
The last condition is equivalent to $i'\ioe i$.  

Property \eqref{TO} follows 
immediately from the construction of the map $\sigma$. 
\epf

\subsubsection{}
The diagrammatic technique for analyzing the term \eqref{jj}
can be formalized as follows 
\begin{equation}
  \label{bMJ}
  \bM(w_1, \ldots,w_m) = \sum_{S\subset \{1,\ldots,m\}} \sum_{\cI} 
\Contr(S,\cI)\,,
\end{equation}
where $\Contr(S,\cI)$ denotes the contribution of the 
interaction diagram $\cI$ to the transformation of term 
indexed by $S$ in \eqref{bM}.

Theorem \ref{tM} is proven by finding a large cancellation 
in the expansion \eqref{bMJ}. We will call an interaction 
diagram $\cI$ \emph{finite} if no interaction with 
infinity occur. The terms in \eqref{bMJ} with $S=\emptyset$
and finite $\cI$ will be called the \emph{principal terms}. 
Theorem \ref{tM} follows from Lemmas \ref{xsxq} and \ref{ysyq} below.

\begin{lm}\label{xsxq}
The minimal $t$ power in the sum of the principal terms is: 
$$t^m \ \sum_{\pi} 
\left(\timeprod_{k=1,\dots,\ell(\pi)} \bA^\vee (w_{\pi_k}) 
\right)\, e^{\alpha_1} \,$$
where the sum is over all partitions $\pi$ of $\{1,\ldots,m\}$.
\end{lm}

\bpf
Let $\cI$ be a finite interaction diagram and let $k$ a
vertex without outgoing edges (that is, a maximal 
element of the associated partial ordering $\ioe$). 
Let 
\begin{equation}
  \label{K}
  K=\{k_1 < k_2 < \dots < k_r=k\} \,,
\end{equation}
be the vertices of $\cI(k)$. 
Since $k<\infty$,
the procedure of Section \ref{alg} requires a selection of a nonminimal
power of $t$ from the operator \eqref{bbA}. Using Lemma \ref{sspp}
we can sum the first nonminimal powers of $t$ over all diagrams
$\cI(k)$ with the vertex set $K$. This yields
\begin{equation}
  \label{t1}
 t^r\,   
\bA^\vee(w_{k_1},\dots,w_{k_r}) \, ,
\end{equation}
completing the proof of the Lemma.
\epf

\begin{lm}\label{ysyq}
The sum of the non-principal terms is 0.
\end{lm}

\bpf
We explain the cancellation of non-principal terms
corresponding to an interaction diagrams $\cI$ 
with exactly one infinite interaction,
\begin{equation}
\label{inff}
k\to \infty.
\end{equation}
The derivation of the general case of the cancellation is
identical.

As before, let \eqref{K} denote the vertex set of $\cI(k)$. 
In contrast to the case considered in Lemma \ref{xsxq},
we now want to extract the \emph{minimal} nonzero power
of $t$ from the operator \eqref{bbA}. Summed over all 
$\cI(k)$ with vertex set $K$, this gives 
\begin{equation}
  \label{t0}
 t^{r-1} \,  \bT(w_{k_1},\dots,w_{k_r})\, 
\bA^0(w_{k_1}+\dots+w_{k_r})  \,. 
\end{equation} 
After commuting past $e^{\alpha_1}$, the 
operator \eqref{t0} is transformed to 
\begin{equation}
 t^{r-1} \,  \bT(w_{k_1},\dots,w_{k_r})\, 
\cE(w_{k_1}+\dots+w_{k_r})  \,.
\end{equation} 

This contribution of the set $K$ is canceled by terms in  
the operator $\bE$  occurring in the definition of $\bM$. 
Consider the summand of formula \eqref{bM} indexed by 
$$
S'= S \cup K.
$$
Consider the associated summand of $\bE(w_S,-t)$ indexed
by 
$$
\pi' = \pi \sqcup K.
$$
It produces the identical result, the only difference 
being that the prefactor of the $K$ contribution,
$$
t^{|K|-1} \, (-1)^{|S|},
$$
is replaced by the factor
$$
(-t)^{|K|-1} \, (-1)^{|S'|}\,,
$$
resulting in a cancellation.
\epf

\label{bzbb}

The following result is
an immediate consequence of
equation \eqref{nloc} and
Theorem \ref{tM}.
Let %
\begin{equation}
  \label{bM'}
  \bM'(w_1,\dots,w_n) = \sum_{\pi} 
\left(\timeprod_{k=1,\dots,\ell(\pi)} \bA^\vee (w_{\pi_k}) 
\right)\, e^{\alpha_1}\,.
 \end{equation} 

\begin{cor}
  \begin{equation}
    \label{neqO}
   \bG'(z_1,\dots,z_n,w_1,\dots,w_m|\nu) = 
\lang \left.
\prod \bA^0(z_i) \, \bM'(w_1,\dots,w_m) \right|
\nu \rang\,. 
\end{equation}
\end{cor}

\subsection{Analysis of the operator formula}\label{sao}

\subsubsection{}
We study here the operator formula \eqref{neqO} for $\bG'$. 
The operator $\bM'$ is constructed from $\bA^\vee$-operators which
in turn are defined in terms of 
 $\bA^0$-operators  and $\bA^1$-operators.
Individual relative invariants of the cap are obtained
by extracting specific powers of the parameters 
$z_i$ and $w_j$.
Therefore, the relative invariants involve the operators
$\bA^0_{i}$ and $\bA^1_{i}$ defined by 
\begin{equation}
\bA^k(z) = \sum_{i\in \Z} \bA^k_i \, z^{i+1} \,, \quad 
k=0,1 \,. 
\label{cbA}
\end{equation}

We will prove a rule for the removal of the $\bA^1$-operators
from formula \eqref{neqO}. The result
expresses the relative invariants of the cap
purely in terms of $\bA^0$-operators. 
By the GW/H correspondence of \cite{OP2}, the vacuum expectations
of $\bA^0$-operators coincide with the stationary invariants of the cap.
We will see in Section \ref{sV} that the rule
for removal of the $\bA^1$-operators 
is equivalent to Virasoro constraints for the cap.

\subsubsection{}

Our strategy for removing the $\bA^1$-operators from 
formula \eqref{neqO} is the following. We commute
all the $\bA^1$-operators, starting from the 
leftmost operator, to the left until they reach the 
vacuum, that is, until they reach the ``$\langle$'' symbol.
Along the way, we will encounter commutators
of $\bA^0$-operators with $\bA^1$-operators.
As we shall see in Lemma \ref{p01}, the commutators 
produce new $\bA^0$-operators. The commutators are related
to the Virasoro reactions of type (iii).

Once a $\bA^1$-operator has reached the vacuum $\vac$,
the $\bA^1$-operator is exchanged for a sum of
$\bA^0$-operators.
Since the vectors
\begin{equation}
\left(\prod \bA^0_{m_i}\right)^* \, \vac\,, \quad 
m_1 \ge m_2 \ge \dots \ge 0 \,, 
\label{A0*}
\end{equation}
span a basis of  $\LV_0$, 
we can express the vector $\left(\bA^1_k\right)^*\, \vac$
as a linear combination of the vectors \eqref{A0*}.
The linear combination, described
in  
Proposition \ref{tV}, is related  to the Virasoro reactions
of type (iv) and (v).

\subsubsection{}

We first study the 
commutation relation between the operators 
$\bA^k(z)$, or, equivalently, their coefficients 
\eqref{cbA}.

Since the operators $\bA^0(z)$ commute for 
different values of $z$,  
$$
\left[\bA^0_k,\bA^0_l\right] = 0 \,.
$$
The commutation of the operators $\bA^0(z)$
can be seen, for example, from equation \eqref{AaE}.

The commutation relation between $\bA^0(z)$ and 
$\bA^1(w)$ can be obtained by extracting the 
$t^1$ term from equation \eqref{bAw}.

\begin{lm}\label{p01} We have
\begin{equation}
    \label{c01}
 \left[\bA^0(z), \bA^1(w)\right] = 
zw\, \bA^0(z+w) \,,
\end{equation}
or, equivalently, 
\begin{equation}
    \label{c01_}
 \left[\bA^0_k, \bA^1_l\right] = 
\binom{k+l}{k} \, \bA^0_{k+l-1} \,, \quad k\ge 0 \,.
\end{equation}
\end{lm}

While $k$ is nonnegative in
equation \eqref{c01_}, 
$k+l-1$ may be negative. From the definitions, 
$$
\bA^0_k  = 
\begin{cases}
1\,, & k=-2 \,, \\
0\,, & k=-1,-3,-4,\dots \, ,
\end{cases}
$$
see \cite{OP2}, \cite{OP3} for a discussion.
Since 
\begin{equation}
\binom{-1}{k} = (-1)^k \,,
\label{bin-1}
\end{equation}
we find,
$$
\left[\bA^0_k, \bA^1_{-k-1}\right] = (-1)^k \,, \quad k\ge 0 \,.
$$
In equation \eqref{bin-1} and below, we view 
$$
\binom{l}{k} = \frac{l(l-1)\dots (l-k+1)}{k!}
$$
as a polynomial of degree $k$ in $l$ which can be 
evaluated at any $l$ and differentiated 
with respect to $l$. 

For the proof of Lemma \ref{lm410} we will 
need the following auxiliary result, namely
the commutation relation between $\bA^0(z)$ and 
$\bA^2(w)$. It can be obtained by extracting the 
$t^2$ term from equation \eqref{bAw}.

\begin{lm} We have
\begin{equation}
    \label{c02}
 \left[\bA^0(z), \bA^2(w)\right] = 
\frac{zw}{1+\frac{z}{w}} \, \bA^1 (z+w)
+
zw^2\, \ln\left(1+\frac zw\right) \, \bA^0(z+w)
 \,,
\end{equation}
or, equivalently,  
\begin{equation}
    \label{c02_}
 \left[\bA^0_k, \bA^2_l\right] = 
\binom{k+l-1}{k} \, \bA^1_{k+l-1} +
\left[\frac{\partial}{\partial l} \binom{k+l-1}{k} \right]\, 
\bA^0_{k+l-2}\, ,  \, \, \, k\geq 0. 
\end{equation}
\end{lm}

Equation \eqref{c02} is obtained from equation \eqref{c02_} by
rewriting the coefficient of the second term,
\begin{eqnarray*}
[z^{k+1} \, w^{l+1}] \, zw^2\, \ln\left(1+\frac zw\right) \, (z+w)^{k+l-1} &= & 
[z^k] \, \ln(1+z) \, (1+z)^{k+l-1} \\ 
& = &\frac{\partial}{\partial l} \, [z^k] \, (1+z)^{k+l-1} \\
& = &\frac{\partial}{\partial l} \binom{k+l-1}{k}\, . \\
\end{eqnarray*}

Finally, we determine the commutation relations
between $\bA^1_k$ and $\bA^1_l$. The commutator 
of $\bA(z)$ with $\bA(w)$, computed in \cite{OP3}, is 
linear the $t$. Hence, 
$$
[t^2] \, \left[ \bA(z), \bA(w) \right] = 0 \,.
$$
The following result is an immediate consequence.

\begin{lm}\label{lm410} We have
\begin{equation}
    \label{c11_}
 \left[\bA^1_k, \bA^1_l\right] = 
\left[\bA^0_l, \bA^2_k\right]
 - \left[\bA^0_k, \bA^2_l\right]\,.
\end{equation}
\end{lm}

\noindent In particular, for $k,l\ge 0$, we find,
$$
\left[\frac{\bA^1_k}{k!}, \frac{\bA^1_l}{l!}\right] =
(k-l) \, \frac{\bA^1_{k+l-1}}{(k+l-1)!} + 
\big(\dots \big) \, \bA^0_{k+l-2} \,.
$$ 
The first term on the right
is
the commutation relation of the Virasoro  
subalgebra $\mathcal{V}$ of the holomorphic vector fields on the 
line (with shifted index).

\subsubsection{}

When 
an $\bA^1$-operator reaches the vacuum 
  ``$\langle$'', we can exchange the $\bA^1$-operator
for $\bA^0$-operators by the following result.

\begin{pr}\label{tV}
We have
\begin{multline}
  \label{1V}
  \lang \bA^1_k \, \dots \rang =  
- \left(\sum_{j=1}^{k}\frac{1}{j}\right) \, \lang \bA^0_{k-1} \, \dots\rang\\
 +
\frac1{2}\,\sum_{i=0}^{k-3} \frac{(i+1)!\,(k-i-2)!}{k!}  
\, \lang \bA^0_i \, \bA^0_{k-i-3} \, \dots \rang \,,
\end{multline}
where the dots stand for arbitrary $\bA^0$ and $\bA^1$-operators. 
\end{pr}

\begin{proof}
For $k\le 0$, 
the vanishing of both sides of equation \eqref{1V}
is obvious.
For $k=1,2,3$, we can check equation \eqref{1V} by
an explicit calculation.

Assume equation \eqref{1V} holds for $k\geq 3$. Then,
using the commutation relation 
$$
\left[\bA^1_2, \bA^1_k\right]  = 
\frac{(k+1)(2-k)}{2} \, \bA^1_{k+1}  + 
\left[(k+1)\sum_{j=2}^{k+1}\frac 1j - k - \frac12\right] \, \bA^0_{k}  \,,
$$
we verify equation \eqref{1V} holds for $k+1$.
\end{proof}

 Lemma \ref{p01} and Proposition \ref{tV} together provide a
rule for the systematic removal of $\bA^1$-operators from
formula \eqref{neqO} for $\bG'$.

\subsection{Virasoro constraints}\label{sV}

\subsubsection{}

To complete the proof of the Virasoro constraints for
the cap, we must prove the formula,
\begin{equation}
    \label{n3}
   \bG'(z_1,\dots,z_n,w_1,\dots,w_m|\nu) = 
\lang \left.
\prod \bA^0(z_i) \, \bM'(w_1,\dots,w_m) \right|
\nu \rang\, ,
\end{equation}
implies the Virasoro constraints for
the depth $m$ theory.
The Virasoro constraints provide rules for 
expressing the depth $m$ theory in
terms of the stationary theory.
We will prove our rule for the removal of
the $\bA^1$-operators from formula \eqref{n3}
yields precisely the {\em same} reduction of the
depth $m$ theory to the stationary theory. 
Therefore,
the unique depth $m$ theory determined by the
Virasoro constraints from the stationary theory 
equals
the depth $m$ theory determined by formula \eqref{n3}.

\subsubsection{}

Consider the reduction of the following relative Gromov-Witten 
invariant,
\begin{equation}
  \label{egw}
  \lang \tau_{l_1}(\omega) \, \tau_{l_2}(\omega) \dots
\tau_{k_1}(1) \, \tau_{k_2}(1) \, \tau_{k_3}(1)
\dots | \nu \,\rang \, , 
\end{equation}
via the Virasoro constraints.
We apply the Virasoro reactions for the removal
of the insertions $\tau_k(1)$ successively 
from left to right starting with $\tau_{k_1}(1)$. 
Of course, the ordering 
in the bracket \eqref{egw} is immaterial. However, since we will interpret
the Virasoro reactions  
in the noncommutative operator formalism of $\LV$, the ordering
will play a crucial role.

When an
insertion $\tau_{k_i}(1)$ decays by interaction with 
another insertion $\tau_r(\gamma)$, the result is placed 
in the location in the bracket 
determined by $\tau_r(\gamma)$. When an insertion 
$\tau_{k_i}(1)$ decays alone, via Virasoro reactions of type (iv) or (v),
the result replaces $\tau_{k_i}(1)$ in the bracket. 

Virasoro reactions of type (iii) in the reduction of
the relative invariant \eqref{egw} can be separated into two subtypes.
A type (iii.a) reaction takes place
when an insertion $\tau_{k_i}(1)$ interacts with a insertion $\tau_r(\omega)$
occurring further {\em right} in the bracket.
A type (iii.b) reaction takes place 
when an insertion $\tau_{k_i}(1)$ interacts with a insertion $\tau_r(\omega)$
occurring further {\em left} in the bracket.

For our analysis, we separate the Virasoro reactions in
the reduction of the relative invariant \eqref{egw}
into two classes:
\begin{enumerate}
\item[I.]
Virasoro reactions of type (i), (ii), and (iii.a),
\item[II.]
Virasoro reactions of type (iii.b), (iv), and (v).
\end{enumerate}
We will first analyze 
the class I interactions. The
 class II interactions will be considered afterwards. We will 
see the total result of all class I interactions
is incorporated in formula \eqref{neqO}.
The class II interactions 
corresponds to the rules for removing the $\bA^1$-operators
derived in Section
\ref{sao}.

\subsubsection{} \label{yo1}

We study here the total effect of all
 class I reactions on the bracket \eqref{egw}.

As a 
first step, we analyze the effects of  
 type (i) reactions. 
Consider the following insertions,
\begin{equation} \label{innt}
\lang \dots \tau_{r_1}(1) \dots \tau_{r_2}(1) \dots 
\tau_{r_s}(1) \dots | \nu \rang,
\end{equation}
in the bracket \eqref{egw}.
The position determined by the last insertion $\tau_{r_s}(1)$
will be called the {\em last} position. A  
 sequence of type (i) reactions which cluster the insertions
\eqref{innt},
the result will be proportional to the bracket
\begin{equation} \label{zqq}
\lang \dots \tau_{\, \sum r_i - s +1}(1) \dots | \nu \rang \,,
\end{equation}
with the new insertion placed in the last position.

\begin{lm} \label{yyjj}
After summation over all possible
sequences of type (i) reactions which cluster the insertions, 
\begin{multline}\label{ri}
\lang \dots \tau_{r_1}(1) \dots \tau_{r_2}(1) \dots 
\tau_{r_s}(1) \dots| \nu \rang  \to \\
  \binom{\sum r_i -1}{r_1,\dots,r_{s-1},r_s-1} \, \lang \dots \tau_{\,
    \sum r_i - s +1}(1) \dots | \nu\rang \,.
\end{multline}
\end{lm}

\bpf The Lemma is easily obtained from
the elementary identity,
\begin{multline*}
  \binom{\sum r_i -1}{r_1,\dots,r_{s-1},r_s-1} =
\binom{\sum r_i -2}{r_2,\dots,r_{s-1},r_1+r_s-2} \,
\binom{r_1+r_s-1}{r_1} + \\
\sum_{j=2}^{s-1} 
\binom{\sum r_i -2}{r_2,\dots,r_1+r_j-1,r_{s-1},r_s-1} \,
\binom{r_1+r_j-1}{r_1}\, ,
\end{multline*}
by induction. \epf

\vspace{10pt}

For fixed values $r_1,\dots,r_{s-1}$, the coefficient in formula \eqref{ri} and
the coefficients of 
the individual type (i) reactions 
are polynomials in $r_s$. Hence, we
may differentiate
the coefficient of formula \eqref{ri}
with respect to the variable $r_s$,
\begin{equation}
\label{dssx}
\frac{\partial}{\partial r_s} \ 
\binom{\sum r_i -1}{r_1,\dots,r_{s-1},r_s-1}.
\end{equation}
The differentiation \eqref{dssx}
provides a summation of over
all Virasoro reactions involving the insertions \eqref{zqq}
of the following form: a sequence of type (i) interactions, followed
by a single type (ii) interaction {\em involving the last position},
followed by a sequence of type (iii.a) interactions.
 
The coefficient in formula \eqref{ri} represents a
sum over all type (i) interactions. 
The differentiation \eqref{dssx} can be
analyzed by considering the individual reactions.
From each sequence of type (i) reactions
involving the insertions \eqref{zqq}, 
the differentiation singles out a particular one {\em involving the last position}
 and
transforms the interaction coefficient by:
$$
\binom{a+b-1}{a} \mapsto \frac{\partial}{\partial b} \,
\binom{a+b-1}{a} = \binom{a+b-1}{a}
\left[\frac1{b}+\dots+\frac1{a+b-1}\right] \,.
$$
where $a$ and $b$ stand for disjoint sums of the
numbers $r_i$ with $r_s$ occurring in $b$.

Lemma \ref{yyjj} and our interpretation of
the differentiation \eqref{dssx} yield a proof of the
following result.

\begin{pr}\label{clI}
After summation over all sequences of class I interactions
which cluster the 
insertions 
$$
\lang \dots \tau_{r_1}(1) \dots \tau_{r_2}(1) \dots 
\tau_{r_s}(1) \dots | \nu \rang, 
$$
we obtain
\begin{multline}\label{rr}
\binom{\sum r_i -1}{r_1,\dots,r_{s-1},r_s-1} \, \lang \dots \tau_{\,
    \sum r_i - s +1}(1) \dots | \nu \rang  + \\
\frac{\partial}{\partial r_s} \, \binom{\sum r_i -1}{r_1,\dots,r_{s-1},r_s-1}
\, \lang \dots \tau_{\,
    \sum r_i - s}(\omega) \dots |\nu \rang \,.
\end{multline}
\end{pr}

\subsubsection{} \label{yo2}

Proposition \ref{clI} has a nice interpretation in terms 
of generating functions. 
First, we have 
\begin{multline*}
  \sum_{r_1,\dots,r_s\ge 0} \binom{\sum r_i
    -1}{r_1,\dots,r_{s-1},r_s-1} \, \tau_{\, \sum r_i - s +1}(1) \, \prod
  w_i^{r_i+1} = \\
\bT(w_1,\dots,w_s)\, 
\frac{w_s}{w_1+\dots+w_s} \,
\sum_{k\ge 0} \tau_k(1) \, \left(\sum w_i\right)^{k+1}\,,
\end{multline*}
which precisely matches the form of the 
$\bA^1$-term in $\bA^\vee$. 
Second, we have 
\begin{align*}
  \frac{\partial}{\partial r_s} \, &\binom{\sum r_i -1}
{r_1,\dots,r_{s-1},r_s-1} \\
&=
\frac{\partial}{\partial r_s}\,
\left[ w_1^{r_1} \, w_2^{r_2} \dots w_{s-1}^{r_{s-1}}\right] \,
\left(\sum_{i=1}^{s-1} w_i + 1 \right)^{\sum_1^s r_i -1} \\
&= 
\left[ w_1^{r_1} \, w_2^{r_2} \dots w_{s-1}^{r_{s-1}}\right] \,
\ln\left(1+\sum_{i=1}^{s-1} w_i\right) \, 
\left(\sum_{i=1}^{s-1} w_i + 1 \right)^{\sum_1^s r_i -1} \\
&=\left[ w_1^{r_1} \, w_2^{r_2} \dots w_s^{r_s-1}\right] \, 
\ln\left(1+\frac{\sum_{1}^{s-1} w_i}{w_s}\right) \, 
\left(\sum_{i=1}^{s} w_i\right)^{\sum_1^s r_i -1} \,.
\end{align*}
Therefore,
\begin{multline*}
  \sum_{r_1,\dots,r_s\ge 0} \frac{\partial}{\partial r_s}\,
\binom{\sum r_i
    -1}{r_1,\dots,r_{s-1},r_s-1} \, \tau_{\, \sum r_i - s}(\omega) \, \prod
  w_i^{r_i+1} = \\
 \bT(w_1,\dots,w_s)\, w_s\,
\ln\left(1+\frac{\sum_{1}^{s-1} w_i}{w_s}\right)\,
\sum_{k\ge 0} \tau_k(\omega) \, \left(\sum w_i\right)^{k+1}\, ,
\end{multline*}
which precisely matches the form of the $\bA^0$ term of $\bA^\vee$.

\subsubsection{}
Consider the generating function of depth $m$ relative invariants
of the cap,
\begin{multline*}
   \bG'(z_1,\dots,z_n,w_1,\dots,w_m|\nu) = \\
 \sum_{k_i\ge 0} \sum_{l_j\ge 0} \, \prod_{i=1}^n z_i^{k_i+1} 
\prod_{j=1}^m w_j^{l_j+1} \, 
\lang \left. \prod \tau_{k_i}(\omega)  \, \prod \tau_{l_j}(1) 
\right| \nu \rang  \,.
\end{multline*}
We reduce the the relative invariants on the right via
all Virasoro reactions of class I.
By the result of Sections \ref{yo1} and \ref{yo2}, the outcome
exactly equals the operator formula,
\begin{equation}
\label{aall}
\lang \left.
\prod \bA^0(z_i) \, \bM'(w_1,\dots,w_m) \right|
\nu \rang\, ,
\end{equation}
after the following substitutions, 
$$
\tau_k(1) \to \bA^1_k \,, \quad 
\tau_k(\omega) \to \bA^0_k \,.
$$

To complete the Virasoro reduction of the depth $m$ theory to the
stationary theory, we must apply Virasoro reactions (iii.b), (iv), and (v). 
By Lemma \ref{p01} and Proposition \ref{tV}, these Virasoro reactions
exactly correspond to the rules for the removal of the $\bA^1$-operators
from formula \eqref{aall}. 

We have proven the reduction of the depth $m$ theory of the cap to
the stationary theory by the Virasoro constraints {\em equals} the 
reduction of formula \eqref{aall} by our operator methods.
Hence, the Virasoro constraints for the cap are proven. \qed

\subsubsection{}
By Proposition \ref{eee}, the proof of the even Virasoro constraints
for the relative theories of all target curves is complete.
The treatment of the odd classes will be presented in
Sections \ref{threeodd} and \ref{three}.

\pagebreak

\section{Odd classes}

\label{threeodd}

\subsection{Overview}
The even relative Gromov-Witten theory of target curves $X$ is completely
determined
by the GW/H correspondence and 
the even Virasoro constraints. 

The full relative Gromov-Witten theory
of target curves  includes the descendents of  both even and odd 
cohomology classes.
We will prove the full relative theory of target curves
is uniquely determined from the even theory
by the following four  properties:
\begin{enumerate}
\item[(i)] Algebraicity of the virtual class,
\item[(ii)] Degeneration formulas for the relative theory in the
presence of odd cohomology,
 \item[(iii)] Monodromy invariance of the relative theory,
\item[(iv)] Elliptic vanishing relations.
\end{enumerate}
The Virasoro constraints for the full theory are proven by establishing
their compatibility with the above properties (i)-(iv).

\subsection{Elliptic invariants}

Let $E$ be an elliptic target with a relative point $e$. 
Let 
$$\alpha \in H^{1,0}(E,\com),$$ 
$$\beta \in H^{0,1}(E,\com),$$
span a symplectic basis with
$$\int_E \alpha \cup \beta = 1.$$
Consider the 
set of relative elliptic invariants with odd insertions:
\begin{equation}
\label{mmmm}
\lang  \left. \prod_{h\in H} \tau_{o_h}(1) 
\prod_{i\in I} \tau_{n_i}(\alpha)
\prod_{j\in J} \tau_{m_j} (\beta) \, \right| \, \eta \rang^E.
\end{equation}
The above invariant is defined by integration against 
$$[\overline{M}_{g,n}(E,\eta)]^{vir}.$$
Since, by property (i), the virtual class is algebraic,
the invariant \eqref{mmmm} vanishes if $|I|\neq |J|$.
The balance $|I|=|J|$ is the {\em only} consequence of the algebraicity which
will be used.
Since the
bracket is skew-symmetric in the odd insertions, we require $I$ and $J$
to be ordered sets to fix the sign.

\begin{pr} 
\label{pqw}
The full relative theory of target curves is uniquely
determined from the even theory by
the elliptic invariants \eqref{mmmm} and the degeneration property (ii).
\end{pr}

\bpf Consider a relative Gromov-Witten invariant on a target curve $X$ of
genus $g$,
\begin{equation}
\label{kssw}
\lang \prod \tau_{r_i}(\gamma_i) \, | \, \eta^1, \ldots, \eta^m \rang^X.
\end{equation}
If $X$ is rational, the theory is even. 

\begin{figure}[!htbp]
  \begin{center}
    \scalebox{0.8}{\includegraphics{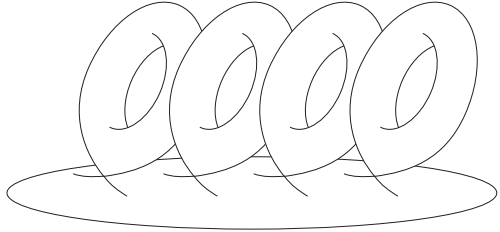}} 
    \caption{A rational curve with 4 elliptic tails}
    \label{f5}
  \end{center}
\end{figure}

If $g>0$, we may degenerate $X$ to a rational curve with  $g$ elliptic
tails, see Fig.~\ref{f5}. 
We may specialize the relative points of $X$ and
the descendents $\tau_r(\omega)$ to the  
rational component.

The degeneration does not alter $H^1(X,\com)$.
The odd cohomology of $X$ can be written in terms of a symplectic basis 
which is the union of the bases of the odd cohomologies of the elliptic components.
The descendents of the odd basis elements then
specialize to the corresponding elliptic factors.

By the degeneration formula, the original invariant \eqref{kssw} is expressed
in terms of the relative invariants of the degenerate components after
all possible distributions of the 
 descendents $\tau_r(1)$.
Hence, the invariant \eqref{kssw} is expressed purely in terms of relative 
invariants of genus
0 from the rational component and invariants 
of  type \eqref{mmmm} from the elliptic components.
\epf

After introducing the monodromy and elliptic vanishing relations in
Sections \ref{mrel} and \ref{evan}, we will
prove the following uniqueness result in Section \ref{u11}.

\begin{pr} \label{bbbbb}
The elliptic invariants \eqref{mmmm} are uniquely determined
from the even theory by properties (i-iv).
\end{pr}

Together, Propositions \ref{pqw} and \ref{bbbbb} show the full relative 
theory of target curves is uniquely determined from the even theory 
by properties (i-iv).

\subsection{Monodromy relations}
\label{mrel}
\subsubsection{}
\label{mmm1}

We will now find relations  for the absolute Gromov-Witten theory
of an elliptic target $E$ obtained from the
monodromy invariance property (iii). 

Using the moduli of elliptic curves, we can find
a monodromy transformation on $H^1(E,\com)$ satisfying:
\begin{equation}
\label{monod}
\alpha \mapsto \alpha, \ \
\beta \mapsto \alpha + \beta.
\end{equation}
In fact, the monodromy group is $SL_2({\mathbb Z})$, but we will only
require the above transformation.

Let $\Psi$ denote the set $\{\psi^0, \psi^1, \psi^2, \dots\}$.
Let $I$ and $J$
be disjoint ordered index sets such that $|I|>0$ and $|I|=|J|$.
Let 
$${\mathbf{n}}: I \rarr {\Psi},\ \ \ i \mapsto \psi^{n_i},$$
$${\mathbf{m}}: J \rarr {\Psi},\ \ \ j \mapsto \psi^{m_j},$$ 
be descendent assignments.

Let $\delta \subset I$ be a subset.
Let $S(\delta)$ denote the set of subsets of 
$I\cup J$ of cardinality $|I|$ containing $\delta$. 
For $D \in S(\delta),$
Let 
$$\tau_{{\mathbf{n}}, {\mathbf{m}}}(D) =
\prod_{i\in I} \tau_{n_i}(\gamma^D_i) 
\prod_{j\in J} \tau_{m_j}(\gamma^D_j).$$
Here, for $\xi \in I\cup J$, 
$$\gamma^D_\xi= \alpha \text{\ or\ } \beta$$ 
if
$\xi \in D$ or $\xi \notin D$ respectively.

The monodromy invariant
monomial insertion,
$$N=\prod_{h\in H} \tau_{o_h}(1) \prod_{h'\in H'}
\tau_{o'_{h'}}(\omega),$$
will be a idle prefactor in the relations below.

\begin{pr} 
\label{jpx}
For every proper subset $\delta \subset I$, the descendent relation,
$$\sum_{D \in S(\delta)}
\lang N  \ \tau_{{\mathbf{n}}, {\mathbf{m}}}(D)   \rang^E_d
=0,$$
holds  for the Gromov-Witten theory
of $E$.
\end{pr}

\bpf The proof is a straightforward application of the monodromy transformation
\eqref{monod}.
Certainly, the invariant 
\begin{equation} \label{vanny}
\lang N\  \prod_{i\in I} \tau_{n_i}(\gamma_i^\delta) \prod_{j\in J}
\tau_{m_j}(\beta) \rang^E_d
\end{equation}
vanishes due to the imbalance of the odd insertions (since $\delta$
is a proper subset of $I$).
The Proposition is obtained by simply applying transformation
\eqref{monod} to the vanishing invariant \eqref{vanny}.
\epf

Let $R_d(N,{\mathbf{n}}, {\mathbf{m}}, \delta)$ denote the monodromy
relation of Proposition \ref{jpx}:
\begin{equation}
\label{wwkkd}
\sum_{D \in S(\delta)}
\lang N\  \tau_{{\mathbf{n}}, {\mathbf{m}}}(D)   \rang^E_d
=0.
\end{equation}

\subsubsection{}

We will require a formal generalization of the monodromy relation
obtained in Section \ref{mmm1}.

Let $\Psi_{\mathbb Q}$ denote the ${\mathbb Q}$-vector space with basis
given by the set $\Psi$.
Let functions ${\mathbf n},{\mathbf m}$ take values in $\Psi_{\mathbb Q}$.
$${\mathbf n}: I \rarr \Psi_{\mathbb Q}, 
\ \ \ {\mathbf m}:J \rarr \Psi_{\mathbb Q}.$$
Instead of simply assigning to each marking in $i\in I$ a descendent
$\psi^{n_i}$, the function ${\mathbf n}$ assigns to $i$ a finite
linear combination,
$$i \mapsto  c^i_{0} \psi^0 + c^i_{1} \psi^1 + c^i_{2} \psi^2 + \dots.$$
Similarly, the function ${\mathbf m}$ assigns to $j$
a finite linear combination,
$$j \mapsto  c ^j_{0} \psi^0 + c^j_{1} \psi^1+ c^j_{2} \psi^2 + \dots.$$

A richer monodromy relation $R_d(N, {\mathbf n}, {\mathbf m}, \delta)$
is defined by setting
$$\tau_{{\mathbf{n}}, {\mathbf{m}}}(D) =
\prod_{i\in I} \Big( \sum_{q\geq 0} c_q^i\tau_{q}(\gamma^D_i) \Big) 
\prod_{j\in J}\Big( \sum_{q\geq 0} c_q^j \tau_{q}(\gamma^D_j) \Big)$$
on the right side of equation \eqref{wwkkd}.
The richer relation $R_d(N, {\mathbf n}, {\mathbf m}, \delta)$ is
proven by expanding and using Proposition \ref{jpx}.

\subsection{Elliptic vanishing relations}
\label{evan}

We present here geometric 
vanishing relations which constrain the absolute Gromov-Witten theory
of $E$.

Let $K$ be an ordered index set satisfying $|K|>0$. 
Let $P$ be a set partition
of $K$ with parts of size at least 2.
Let $P_1, \ldots, P_\ell$ be the parts of $P$.

Let $\overline{M}_{g,S}(E,d)$ be a moduli space of stable maps 
with possibly disconnected domains for which
the marking set $S$ contains $K$.  Let
$$\phi_i: \overline{M}_{g,S}(E,d) \rarr E^{|P_i|}$$
be the product evaluation map determined by the ordered part $P_i$.
Let $${\mathbf l} : K \rarr \Psi, \ \ \ k \mapsto \psi^{l_k},$$ be
a descendent assignment.

The {\em small diagonal} of the $r$-fold product $E^r$ is defined by:
 $$\{ (x,\ldots,x) \ | \ x\in E \} \subset E^r.$$ 
Let $\bigtriangleup_r\in H^*(E^r,\com)$ denote the Poincar\'e dual of the
small diagonal.

The monomial insertion of descendents of the identity,  
$$M=\prod_{h\in H} \tau_{o_h}(1),$$
will be a idle prefactor in the elliptic vanishing relations below.
The descendents $\tau_k(\omega)$ do {\em not} appear in $M$.

\begin{pr}
The elliptic vanishing relation $V_d(M, P, {\mathbf l})$ holds:
$$\int_{[\overline{M}_{g,H+K}(E,d)]^{vir}} \prod_{h\in H} \psi_h^{o_h} \ 
\prod_{k \in K} \psi_k^{l_k} \ 
\prod_{i=1}^{\ell} \phi_i^*(\bigtriangleup_{|P_i|}) \ \ = 0.$$
\end{pr}

While the Proposition is true for all $g$, the vanishing is trivial unless
$g$ is determined from the rest of the data by 
the dimension constraint. 

\vspace{+10pt}
\bpf
The integral is proven to vanish in two steps. First, the virtual
fundamental class is analyzed. Second, the integrand is analyzed.
A similar elliptic vanishing is proven by the same method in  \cite{pandeg}.

The moduli space of maps $\overline{M}_{g, H+K}(E,d)$ is equipped
with an algebraic translation action of the elliptic curve $E$.
There exists an algebraic quotient of this free action:
$$ \overline{M}_{g, H+K}(E,d)/ E =  \text{ev}_\xi^{-1}(0) 
\subset \overline{M}_{g,H+K}(E,d),$$
where $\xi$ is any marking and $0\in E$ is the neutral element. In fact
$\overline{M}_{g, H+K}(E,d)$ is $E$-equivariantly isomorphic
to a product of $E$ by the quotient.

The
virtual fundamental class of $\overline{M}_{g, H+K}(E,d)$
is pulled-back from the quotient $\overline{M}_{g, H+K}(E,d)/ E$. The pull-back 
property is obtained easily from the
construction of the virtual fundamental class. 
Since no $\tau_k(\omega)$ insertions are allowed,
the integrand is also pulled-back from the quotient space. 
Hence, we may use the push-pull formula to conclude the integral vanishes.
\epf 

The elliptic vanishing relations can be expressed in terms of
the absolute Gromov-Witten theory of $E$ via the K\"unneth decompositions
of the classes $\bigtriangleup_r$.
The K\"unneth decompositions of $\bigtriangleup_2$ and $\bigtriangleup_3$
are
$$\bigtriangleup_2 
= 1\otimes \omega + \omega 
\otimes 1 - \alpha \otimes \beta + \beta \otimes \alpha.$$
\begin{eqnarray*}
\bigtriangleup_3 & = &
1\otimes \omega \otimes \omega + \omega 
\otimes 1 \otimes \omega + \omega \otimes \omega \otimes 1  \\
& & -\alpha \otimes \beta \otimes \omega + \beta \otimes 
\alpha \otimes \omega \\
& & - \alpha \otimes \omega \otimes \beta + 
\beta \otimes \omega \otimes \alpha \\
& & - \omega \otimes \alpha \otimes \beta + \omega 
\otimes \beta \otimes \alpha.
\end{eqnarray*}
We will separate the K\"unneth components
of $\bigtriangleup_r$ into two groups,
$$\bigtriangleup_r = \bigtriangleup_r^{even} + \bigtriangleup_r^{odd}.$$
The summand $\bigtriangleup_r^{even}$
consists of $r$ terms in which the identity class
 and $r-1$ copies of $\omega$ are tensored in all distinct orders,
$$\bigtriangleup_r^{even} = 1\otimes \omega^{r-1}+ \ldots + \omega^{r-1} \otimes 1.$$
The summand $\bigtriangleup_r^{odd}$
consists of $2\binom{r}{2}$ terms. For each pair of indices $i<j$
two terms occurs: 
\begin{enumerate}
\item[$\bullet$] $-\alpha$ in the $i$th factor, $\beta$ in the $j$th factor and
$r-2$ copies of $\omega$ in all the other tensor factors,
\item[$\bullet$] $\beta$ in the $i$th factor, $\alpha$ in the $j$th factor and
$r-2$ copies of $\omega$ in all the other tensor factors.
\end{enumerate}
We will be primarily interested in the summand $\bigtriangleup_r^{odd}$.

The simplest example occurs when $|K|=2$ and $P$ has one part. After expanding
$V_d(M,P,{\mathbf l})$ using the K\"unneth decomposition of $\bigtriangleup_2$,
we find
\begin{eqnarray*}
& & \lang M \ \tau_{l_1}(1) \tau_{l_2}(\omega) \rang^E_d  + 
\lang M \ \tau_{l_1}(\omega) \tau_{l_2}(1) \rang^E_d  \\
&-& \lang M \ \tau_{l_1}(\alpha) \tau_{l_2}(\beta) \rang^E_d 
+ \lang M \ \tau_{l_1}(\beta) \tau_{l_2}(\alpha) \rang^E_d =0. 
\end{eqnarray*}
Descendents of the odd cohomology of $E$
appear via the summand 
$\bigtriangleup_2^{odd}$.

The function ${\mathbf l}$ may take more general values
for the elliptic vanishing relation, 
$${\mathbf l}: K \rarr \Psi_{\mathbb Q}, 
\ \ \ k \mapsto c^k_{0} \psi^0 + c^k_{1} \psi^1 + c^k_{2} \psi^2 + \dots.$$
The relation $V_d(M,P,{\mathbf l})$ is well-defined and true
in the richer context.

\subsection{Proof of Proposition \ref{bbbbb}}
\label{u11}

\subsubsection{}

We must determine the 
relative elliptic invariants,
\begin{equation}
\label{mmmmm}
\lang \prod_{h\in H} \tau_{o_h}(1) 
\prod_{i\in I} \tau_{n_i}(\alpha)
\prod_{j\in J} \tau_{m_j} (\beta) \ | \ \eta \rang^E_d,
\end{equation}
from the even theory by properties (i-iv).

By property (i), the invariants vanish unless
$|I|=|J|$.
We will proceed by induction on $|I|$.

If $|I|=0$, then the invariant is even.
We will start with a proof of Proposition \ref{bbbbb}
in case $|I|=1$. The method for $|I|=1$ will be generalized
in Section \ref{u33} to establish the induction step.

\begin{lm} 
\label{bscs}
The elliptic invariants \eqref{mmmmm} where $|I|=1$ are
uniquely determined from the even theory by properties (i)-(iv).
\end{lm}

For the proof of Lemma \ref{bscs},  we will require an auxiliary result 
derived 
from the GW/H correspondence.

Let ${\mathcal P}(d)$ be the set of
partitions of $d$. Let ${\mathbb Q}^{{\mathcal P}(d)}$ denote the linear
space of functions from ${\mathcal P}(d)$ to
${\mathbb Q}$.
Let $$\tilde{\tau}(\omega)= \sum_{q=0}^\infty c_q \tau_q(\omega)$$ be a {\em finite}
linear combination of descendent of $\omega$.
For $\w\geq 0$, define a function on ${\mathcal P}(d)$ by:
$$\gamma_{\w}: {\mathcal P}(d) \rarr {\mathbb Q}, \ \ \
\eta \mapsto \lang  \tilde{\tau}(\omega)^\w\ |\ \eta\rang_d^{\proj^1}.$$
The above bracket is defined by a multilinear expansion of the
insertion $\tilde{\tau}(\omega)^\w$.

\begin{lm} \label{cakeshop}
For $d\geq 0$, there exists a linear combination $\tilde{\tau}(\omega)$ for which
the set of functions,  
$$\{\gamma_{0}, \gamma_{1}, \gamma_{2}, \dots \},$$
spans ${\mathbb Q}^{{\mathcal P}(d)}$.
\end{lm}

\bpf The formula, 
$$\gamma_{\w}(\eta) = \sum_{|\lambda| = d} 
\left( \frac{ \text{dim} \lambda}{d!} \right) ^2  
\left(
\sum_{q=0}^\infty c_q\frac{ {\mathbf {p}}_{q+1}(\lambda)}{(q+1)!} 
\right)^\w
{\mathbf f}_\eta(\lambda),$$
is a direct consequence of the GW/H correspondence \cite{OP2}.

Since ${\mathbf f}_{\eta}(\lambda)$ is proportional to 
the character of the conjugacy class $C_\eta$ in the representation $\lambda$
of the symmetric group $S_d$,
$${\mathbf f}_{\eta}(\lambda) = |C_\eta| \frac{\chi^\lambda_\eta}{\text{dim} \lambda},$$
the functions, 
 $$\eta \mapsto {\mathbf f}_\eta(\lambda),$$
span ${\mathbb Q}^{{\mathcal P}(d)}$
as $\lambda$ varies. 

To prove the Lemma, we require a $\tilde{\tau}(\omega)$ for
 which the functions,
$$
\lambda \mapsto 
\left( \sum_{q=0}^\infty c_q \frac{{\mathbf {p}}_{q+1}(\lambda)}
{(q+1)!}\right)^\w\,,
$$
span
${\mathbb Q}^{{\mathcal P}(d)}$ as $\w$ varies.
By the Vandermonde determinant, we need only find a $\tilde{\tau}(\omega)$ for which
the values 
$$\sum_{q=0}^\infty c_q \frac{{\mathbf {p}}_{q+1}(\lambda)}{(q+1)!}$$
are distinct as $\lambda$ varies in ${\mathcal P}(d)$.

Since $\lambda$ is a partition of $d$, we may write 
$$\lambda= \lambda_1 \geq \lambda_2 \geq \dots \geq \lambda_{d} \geq 0.$$
We may recover $\lambda$ from the set:
\begin{equation}
\label{elem}
\lambda_1 - 1 +\frac{1}{2}, 
\lambda_2-2+\frac{1}{2}, \ldots, \lambda_d-d+\frac{1}{2}.
\end{equation}

On ${\mathcal P}(d)$, the function ${\mathbf p}_1$ is easily evaluated
to yield a nonzero constant. 
By definition, the functions 
 ${\mathbf p}_{q+1}(\lambda)$ are (up to constants) the $q+1$-power sums of the
elements \eqref{elem}.
Since the functions ${\mathbf p}_1, {\mathbf p}_2, \dots$ include all the power sums,
their values separate elements of ${\mathcal P}(d)$.
Since ${\mathcal P}(d)$ is a finite set, we may find a finite
linear combination of elements ${\mathbf p}_1,{\mathbf p}_2, \dots$ which separate the
set.
\epf

We now prove Lemma \ref{bscs}. We will start by proving the invariants
\begin{equation}
\label{mmmmmx}
\lang  
 \tau_{n}(\alpha)
 \tau_{m} (\beta)\ | \ \eta \rang^E_d
\end{equation}
are determined from the even theory by properties (i)-(iv).

Let $d\geq 0$.
Let $\tilde{\psi} = \sum_{q\geq 0} c_q \psi^q$, where
$$\tilde{\tau}(\omega) = \sum_{q\geq 0} c_q \tau_q(\omega)$$
satisfies the conditions of Lemma \ref{cakeshop} for $d$.

Let $\w\geq 0$. Let $K_\w$ be an ordered index set with $\w+2$ elements.
Let $P$ be the set partition of $K_\w$ with one part.
Let the descendent assignment
${\mathbf l}$ on $K_\w$ take the value $\tilde{\psi}$ for all elements of $K_\w$.

Consider the elliptic vanishing relation $V_d(1,P, {\mathbf l})$.
The terms of $V_d(1,P,{\mathbf l})$ which contain odd classes from the
K\"unneth decomposition of $\bigtriangleup_{\w+2}$ are easily seen to equal:
\begin{eqnarray*}
& - & \binom{\w+2}{2} \lang \tilde{\tau}(\omega)^\w \ \tilde{\tau}(\alpha) 
\tilde{\tau}(\beta)\rang_d^E\\
& +& \binom{\w+2}{2} \lang \tilde{\tau}(\omega)^\w \ \tilde{\tau}(\beta) 
\tilde{\tau}(\alpha)\rang_d^E.
\end{eqnarray*}
After an application of the monodromy relation 
$R_d( \tilde{\tau}(\omega)^\w, 
\{\tilde{\psi}\}, \{\tilde{\psi}\}, \emptyset),$ we may rewrite the
odd terms as:
\begin{eqnarray*}
& - 2& \binom{\w+2}{2} \lang \tilde{\tau}(\omega)^\w \ \tilde{\tau}(\alpha) 
\tilde{\tau}(\beta)\rang_d^E.
\end{eqnarray*}
As the remainder of the relation $V_d(1,P,{\mathbf l})$
consists of terms with only even descendent insertions, we 
may conclude the invariants
\begin{equation}
\label{sss}
\lang \tilde{\tau}(\omega)^\w \ \tilde{\tau}(\alpha) 
\tilde{\tau}(\beta)\rang_d^E
\end{equation}
are determined for all $\w$.

We now study the invariants \eqref{sss} via the degeneration formula:
\begin{equation}
\label{degg}
\lang \tilde{\tau}(\omega)^\w \ \tilde{\tau}(\alpha) 
\tilde{\tau}(\beta)\rang_d^E = \sum_{|\eta|=d} 
\lang \tilde{\tau}(\alpha) 
\tilde{\tau}(\beta)\ | \ \eta \rang_d^E
\ \zz(\eta)\
\lang  \eta \ | \ 
\tilde{\tau}(\omega)^\w \rang_d^{\proj^1},
\end{equation}
where $\zz(\eta)= |\text{Aut}(\eta)| \prod_i \eta_i$.
Here, $E$ degenerates to a nodal target $$E \cup \proj^1.$$ 
The $\w$ markings corresponding to the insertions $\tilde{\tau}(\omega)$
specialize to the component $\proj^1$ in the degeneration.

We have seen  the left side of \eqref{degg}
is determined for all $\w$ from the even theory by conditions (i)-(iv). 
The invariants
\begin{equation}
\label{sttp}
\lang  \tilde{\tau}(\alpha) 
\tilde{\tau}(\beta)\ | \eta \rang_d^E
\end{equation}
are then uniquely determined by Lemma \ref{cakeshop}.

Let $L$ be an arbitrary monomial in the descendents of $\omega$,
$$L = \prod_{h'\in H'} \tau_{o'_{h'}}(\omega).$$
By the degeneration formula, the invariants
\begin{equation}
\label{ksk}
\lang L \ \tilde{\tau}(\alpha)\tilde{\tau}(\beta) \rang^E_d
\end{equation}
are determined by the invariants \eqref{sttp} and the
even relative theory of $\proj^1$.

Let $K_\w$ and $P$ be as before.
Let ${\mathbf l}_{f}$ take the value $\psi^n$
on the first element of $K_\w$ and the value $\tilde{\psi}$
on the following elements.
Consider the elliptic vanishing relation $V_d(1,P, {\mathbf l}_f)$.
The terms of $V_d(1,P,{\mathbf l}_f)$ which contain odd classes from the
K\"unneth decomposition are:
\begin{eqnarray*}
& - & \binom{\w+1}{1} \lang \tilde{\tau}(\omega)^\w \ \tau_n(\alpha) 
\tilde{\tau}(\beta)\rang_d^E\\
& +& \binom{\w+1}{1} \lang \tilde{\tau}(\omega)^\w \ \tau_n(\beta) 
\tilde{\tau}(\alpha)\rang_d^E \\
& - & \binom{\w+1}{2} \lang \tilde{\tau}(\omega)^{\w-1} \tau_n(\omega) \ 
\tilde{\tau}(\alpha)\tilde{\tau}(\beta)\rang_d^E \\
& +& \binom{\w+1}{2} \lang \tilde{\tau}(\omega)^{\w-1} \tau_n(\omega) \ 
\tilde{\tau}(\beta)\tilde{\tau}(\alpha)\rang^E_d.
\end{eqnarray*}

By the
determination \eqref{ksk}, only the first two terms need be analyzed.
After application of the monodromy relation 
$$R_d( \tilde{\tau}(\omega)^\w, 
\{\ \psi^n \}, \{\tilde{\psi}\}, \{1\}),$$ the first two odd terms 
equal:
\begin{eqnarray*}
& - 2& \binom{w+1}{1} \lang \tilde{\tau}(\omega)^\w \ \tau_n(\alpha) 
\tilde{\tau}(\beta)\rang_d^E.
\end{eqnarray*}
As the remainder of the relation $V_d(1,P,{\mathbf l}_f)$
consists of even terms, we may conclude the invariants
\begin{equation*}
\lang \tilde{\tau}(\omega)^\w \ {\tau}_n(\alpha) 
\tilde{\tau}(\beta)\rang_d^E
\end{equation*}
are determined for all $\w$.

Now, by degeneration and Lemma \ref{cakeshop}, as before, we find the
invariants,
$$\lang {\tau}_n(\alpha) 
\tilde{\tau}(\beta) \ | \ \eta \rang_d^E, \ \ \
\lang  L \ {\tau}_n(\alpha) 
\tilde{\tau}(\beta)\rang_d^E,$$
are determined.

Similarly, by studying the  elliptic vanishing relation 
$V_d(1,P, {\mathbf l}_l)$ for the
function ${\mathbf l}_l$ which takes the value  $\psi^m$ on the last element of
$K_\w$ and $\tilde{\psi}$ on the preceding elements, we find the
invariants,
$$\lang   \tilde{\tau}(\alpha) 
{\tau}_m(\beta)\ | \ \eta \rang_d^E, \ \ \
\lang  L \ \tilde{\tau}(\alpha) 
{\tau}_m(\beta)\rang_d^E,$$
are determined.

Finally, we study the elliptic vanishing relations $V_d(1,P, {\mathbf l}_{fl})$
where the function
${\mathbf l}_{fl}$ takes the value $\psi^n, \psi^m, \tilde{\psi}$ on the first,
last, and remaining elements of $K_\w$ respectively.
We then find the invariants,
$$\lang   {\tau}_n(\alpha) 
{\tau}_m(\beta)\ | \eta \rang_d^E,$$ 
are determined.

To conclude the proof of Proposition \ref{bscs}, we must show the invariants
$$\lang \eta \, | \, M \ \tau_n(\alpha) \tau_m(\beta)\rang_d^E$$
are determined for every monomial $M= \prod_{h\in H} \tau_{o_h}(1)$.

We proceed by induction on the degree of $M$. The degree 0 case 
has already been established. If $\text{deg}(M)>0$, we observe that $M$
is a spectator in both the monodromy and elliptic vanishing relations.
Hence, we may repeat the above argument 
based upon the elliptic vanishing relations 
$$V_d(M,P, {\mathbf l}), \, V_d(M,P,{\mathbf l}_f), \, V_d(M,P,{\mathbf l}_l), \,
V_d(M,P, {\mathbf l}_{fl}),$$
where the definitions of the functions ${\mathbf l}$, ${\mathbf l}_f$,
${\mathbf l}_l$, and ${\mathbf l}_{fl}$ on $K_\w$ are unchanged.

The only difference occurs in the degeneration formulas. Here, we must
sum over all possible distributions of $M$ over the degenerate components.
However, if any factors of $M$ are distributed to the component $\proj^1$, 
the
resulting $\tau_k(1)$ monomial on the component $E$ will have strictly lower
degree. Consequently, the terms in which factors of $M$ are distributed to
$\proj^1$ are inductively determined. Hence, we need only consider
terms in the degeneration formula for which the entire monomial $M$
remains on the component $E$. Then, the induction step proceeds exactly as
the degree 0 case.
\epf

\subsubsection{}

\label{u22}

We will derive consequences of the
monodromy and elliptic vanishing relations needed for the
$|I|$ induction in the proof of Proposition \ref{bbbbb}.

Let $I,J$ be disjoint ordered index sets satisfying $|I|=|J|$.
Let ${\mathbf n}, {\mathbf m}$ be functions,
$${\mathbf n}: I \rarr \Psi,\ \ {\mathbf m}:J \rarr \Psi.$$
Let $K= I \cup J$. We order $K$ by placing $I$ before $J$. Let
$${\mathbf l}: K \rarr \Psi$$
be determined by ${\mathbf n}, {\mathbf m}$.

We will consider two types of relations.
Let $\sigma: I \rarr J$ be a bijection. 
Let $P_\sigma$ be the set partition of $K$ into
doublets given by $\{i, \sigma(i)\}$.
First, as $\sigma$ varies, we find $|I|!$ relations,
\begin{equation}
\label{ooo}
V_d(M,P_\sigma, {\mathbf l}),
\end{equation}
where $M= \prod_{h\in H} \tau_{o_h}(1)$ is a fixed monomial.
Second, we have all the monodromy relations,
\begin{equation}
\label{ttttt}
R_d(M, {\mathbf n}, {\mathbf m}, \delta),
\end{equation}
for proper subsets $\delta \subset I$.

\begin{lm} \label{jjj}
The relations \eqref{ooo} and \eqref{ttttt} determine
the invariant,
$$\lang M \ \prod_{i\in I} \tau_{n_i}(\alpha) \prod_{j\in J} \tau_{m_j}(\beta)
\rang^E_d,$$
in terms of degree $d$ invariants of $E$ with strictly fewer odd insertions.
\end{lm}

\bpf
Let $\delta \subset I$ be a subset, and 
let $S(\delta) $ be the set of subsets 
$I\cup J$ of cardinality $|I|$ containing $\delta$.
For $D \in S(\delta),$
let 
$$\tau_{{\mathbf{n}}, {\mathbf{m}}}(D) =
\prod_{i\in I} \tau_{n_i}(\gamma^D_i) 
\prod_{j\in J} \tau_{m_j}(\gamma^D_j),$$
following the notation of Section \ref{mrel}.
%
%
Let $S^*(\delta) \subset S(\delta)$ denote the
set of subsets $D$ satisfying $D \cap I = \delta$.

Since we are only interested in invariants with $|I|+|J|$ odd insertions,
we need only analyze the odd splittings of the $|I|$ distinct
K\"unneth decompositions
in the relation $V_d(M, P_\sigma, {\mathbf l})$. Since $I$ and $J$ are
ordered, the function $\sigma$ is canonically an element of the
symmetric group and therefore has a sign.
We easily compute 
the sum of the terms of 
\begin{equation} \label{pws}
 \sum_{\sigma} (-1)^{\binom{|I|}{2}}  \text{sign}(\sigma) \ V_d(M,P_\sigma,{\mathbf l})
\end{equation} 
with $|I|+|J|$ odd parts equals:
\begin{equation}
\label{bbgg}
\sum_{\delta \subset I}\ 
\sum_{D \in S^*(\delta)} \ 
(-1)^{|I|-|\delta|}|\delta|! (|I|-|\delta|)!\ 
\lang M \ \prod_{i\in I} \tau_{n_i}(\gamma^D_i) \prod_{j\in J}
\tau_{m_j}(\gamma^D_j) \rang^E_d.
\end{equation}

The invariant
$\lang M \ \prod_{i\in I} \tau_{n_i}(\alpha) \prod_{j\in J} \tau_{m_j}(\beta)
\rang^E_d$ occurs in the $\delta= I$ summand of \eqref{bbgg} with
coefficient $|I|!$ as the invariant appears exactly once in each
summand of \eqref{pws}. Let $V$ denote the sum \eqref{bbgg}.

For every $\ell <|I|$, let $R(\ell)$ denote
the monodromy relation sum,
$$\sum_{|\delta| =\ell} R_d(M, {\mathbf n}, {\mathbf m}, \delta).$$
We may expand $R(\ell)$ as:
\begin{equation}
\label{mmgg}
\sum_{|\delta| \geq \ell}\ 
\sum_{D \in S^*(\delta)} \ 
\binom{|\delta|}{\ell}\ 
\lang M \ \prod_{i\in I} \tau_{n_i}(\gamma^D_i) \prod_{j\in J}
\tau_{m_j}(\gamma^D_j) \rang^E_d\ = \ 0.
\end{equation}
 
Using the relations $R(0), \dots, R(|I|-1)$, we can uniquely eliminate
all terms of $V$ except for the $\delta= I$ term,
\begin{equation}\label{sqa}
\lang M \ \prod_{i\in I} \tau_{n_i}(\alpha) \prod_{j\in J} \tau_{m_j}(\beta)
\rang^E_d,
\end{equation}
which we hope to determine. If the coefficient of the term \eqref{sqa} is
not 0 after elimination, then the Lemma is proven.

We abstract the linear algebra arising in the above elimination.
Let ${\mathbb Q}^{|I|+1}$ be
a vector space with basis $e_0, e_1, \dots, e_{|I|}$.
Here, $e_k$ corresponds to the sum,
\begin{equation*}
\sum_{|\delta| = k}\ 
\sum_{D \in S^*(\delta)} \  
\lang M \ \prod_{i\in I} \tau_{n_i}(\gamma^D_i) \prod_{j\in J}
\tau_{m_j}(\gamma^D_j) \rang^E_d .
\end{equation*}
Then,
$V$ is the vector,
$$V = \sum_{k =0}^{|I|} (-1)^{|I|-k} k! (|I|-k)!  \ e_k.$$
For $0 \leq \ell \leq |I|$, let
$$R(\ell)= \sum_{k \geq \ell}  \binom{k}{\ell} \ e_k.$$
For $\ell<|I|$, the vector $R(\ell)$ is the corresponding the
monodromy relation.

The vectors $R(0), \dots, R(|I|)$ span a basis of ${\mathbb Q}^{|I|+1}$.
Hence,
$$V = \sum_{\ell=0}^{|I|} c_\ell R(\ell),$$
for unique coefficients $c_\ell$.
The coefficient of $e_{|I|}$ obtained after the canonical elimination of
$V$ by the vectors $R(0), \dots, R(|I|-1)$ is simply $c_{|I|}$.

The column vectors $R(\ell)$ determine an $|I| \times |I|$ lower
unitriangular matrix $R$
with coefficients $$R_{ab}=\binom{a}{b}.$$
It is well known that $R^{-1}$ has coefficients
$$(R^{-1})_{ab} = (-1)^{a+b} \binom{a}{b}.$$
Since the column vector $(c_0, \dots, c_{|I|})$ is obtained
by the product of $R^{-1}$ with the column vector $V$, we find
$$c_{|I|} = \sum_{k=0}^{|I|}  (-1)^{|I|+k} \binom{|I|}{k}(-1)^{|I|-k} k! (|I|-k)!
= (|I|+1)!\, .$$
The proof of the Lemma is complete.
\epf

Lemma \ref{jjj} is valid in case the function ${\mathbf n}$ and
${\mathbf m}$ take more general values,
$${\mathbf n}: I \rarr \Psi_{\mathbb Q},\ \ 
 {\mathbf m}:J \rarr \Psi_{\mathbb Q},$$
since the monodromy and elliptic vanishing relations remain valid.

\subsubsection{}

\label{u33}
Consider the relative elliptic invariant
\begin{equation} \label{yqpd}
\lang \prod_{h\in H} \tau_{o_h}(1) 
\prod_{i\in I} \tau_{n_i}(\alpha)
\prod_{j\in J} \tau_{m_j} (\beta) \ | \ \eta \rang^E_d.
\end{equation}
Assume  such invariants with
strictly fewer odd insertions are determined from the
even theory by properties (i)-(iv).
We now complete the proof of Proposition \ref{bbbbb}
by establishing the induction step.
We will follow the proof
of Lemma \ref{bscs} using a variation 
of Lemma \ref{jjj} for the monodromy and elliptic
vanishing relations.

We will start by assuming the monomial, $$M= \prod_{h\in H} \tau_{o,h}(1),$$ 
is degree 0
and proceed by induction on the degree of $M$.

Let $\w\geq 0$.
Let $W$ be an ordered set disjoint from $I$ and $J$
satisfying $|W|=\w$. Let $K_\w$
be defined by
$$K_\w= I \cup W \cup J,$$
with the given order.
Let $1\in I$ denote the first element.
For each bijection $$\sigma: I \rarr J,$$ let $P_\sigma$ be the
set partition given by the part $\{1\} \cup W \cup \{\sigma(1)\}$ of
order $\w+2$ together
with the doublets $\{i, \sigma(i)\}$ for $i\neq 1$.
Let the function ${\mathbf l}$ take the value $\tilde{\psi}$ on
all elements of $K_\w$.

Let $V$ be the sum of the terms of
\begin{equation} \label{hjh}
\sum_\sigma 
\binom{\w+2}{2}^{-1} 
 (-1)^{\binom{|I|}{2}}  \text{sign}(\sigma) \ V_d(M,P_\sigma,{\mathbf l})
\end{equation} 
with $|I|+|J|$ odd parts.
The inverse binomial prefactor accounts for the multiplicity of choice
in the K\"unneth decomposition absent for the doublets considered
in Lemma \ref{jjj}. We find, $V$ equals
\begin{equation}
\label{bbggg}
\sum_{\delta \subset I}\ 
\sum_{D \in S^*(\delta)} \ 
(-1)^{|I|-|\delta|}|\delta|! (|I|-|\delta|)!\ 
\lang M  \tilde{\tau}(\omega)^\w \prod_{i\in I} \tilde{\tau}(\gamma^D_i) 
\prod_{j\in J}
\tilde{\tau}(\gamma^D_j) \rang^E_d,
\end{equation}
following the notation of the proof of Lemma \ref{jjj}.

Next, the monodromy relations
$R(0), \dots, R(|I|-1)$ are considered with the induced descendent assignments
${\mathbf n}, {\mathbf m}$ and prefactor
$M\tilde{\tau}(\omega)^\w$. Elimination proves all the invariants
$$\lang M  \tilde{\tau}(\omega)^\w \ \prod_{i\in I} \tilde{\tau}(\alpha) 
\prod_{j\in J}
\tilde{\tau}(\beta) \rang^E_d,$$
are inductively determined. The elimination analysis exactly follows
the proof of Lemma \ref{jjj}.

Degeneration, together with 
Lemma \ref{cakeshop} and induction on the degree of $M$, then shows all
the invariants
$$\lang M  \ \prod_{i\in I} \tilde{\tau}(\alpha) 
\ \prod_{j\in J}
\tilde{\tau}(\beta)\ | \eta \rang^E_d,$$
\begin{equation} \label{sominv}
\lang M L \ \prod_{i\in I} \tilde{\tau}(\alpha) 
\ \prod_{j\in J}
\tilde{\tau}(\beta) \rang^E_d,
\end{equation}
are inductively determined. Here, $L= \prod_{h'\in H'} \tau_{o'_{h'}}(\omega)$
is an arbitrary monomial.

We will now repeat the analysis for several different assignment functions.
Let ${\mathbf l}_{f[r]l[s]}$ take the values
$${\mathbf l}_{f[r]l[s]}(\xi)= n_\xi$$
for the first $r$ elements of $I$, and the values  
$${\mathbf l}_{f[r]l[s]}(\xi)= m_\xi$$
for the first $s$ elements of $J$, and the value
$\tilde{\psi}$ for the remaining elements of $K_\w$.
We have already considered the function
${\mathbf l}_{f[0]l[0]}$.

We first repeat the analysis for the assignment function
${\mathbf l}_{f[1]l[0]}$.
Let $V$ be  the sum of the terms of
\begin{equation*} 
\sum_\sigma \binom{\w+1}{1}^{-1} 
 (-1)^{\binom{|I|}{2}}  \text{sign}(\sigma) \ V_d(M,P_\sigma,{\mathbf l}_{f[1]l[0]})
\end{equation*} 
with $|I|+|J|$ odd parts {\em modulo the invariants \eqref{sominv}}.

Consider the K\"unneth decomposition associated
to the first part of $P_\sigma$.
The terms with an odd class  distributed to $1$ 
contribute to $V$ (and are normalized by the prefactor
$\binom{\w+1}{1}^{-1}$ since they occur with multiplicity).
If the odd parts are distributed away from $1$,
then the resulting terms
are of the form \eqref{sominv}.

We may then eliminate $V$ using the relations $R(0), \dots, R(|I|-1)$
with the induced descendent assignments and prefactor $M \tilde{\tau}(\omega)^\w$.
Because of the normalization and the removal of the
invariants \eqref{sominv}, the elimination analysis exactly follows the
proof of 
 Lemma \ref{jjj}. 

By degeneration, Lemma \ref{cakeshop}, and induction on the
degree of $M$, we conclude 
the invariants
$$\lang 
  M  \ \tau_{n_1}(\alpha) \prod_{1\neq i\in I} \tilde{\tau}(\alpha) 
\ \prod_{j\in J}
\tilde{\tau}(\beta)\ |\ \eta \rang^E_d,$$
\begin{equation}\label{sominv2}
\lang M L \ \tau_{n_1}(\alpha)\prod_{1\neq i\in I} \tilde{\tau}(\alpha) 
\ \prod_{j\in J}
\tilde{\tau}(\beta) \rang^E_d,
\end{equation}
are inductively determined for any $n_1$.

Next, we repeat the analysis for the assignment function
${\mathbf l}_{f[0]l[1]}$.
Let $V$ be  the sum of the terms of
\begin{equation*} 
\sum_\sigma C^{-1}_\sigma 
(-1)^{\binom{|I|}{2}}  \text{sign}(\sigma) \ V_d(M,P_\sigma,{\mathbf l}_{f[0]l[1]})
\end{equation*} 
with $|I|+|J|$ odd parts {\em modulo the invariants \eqref{sominv}}.

Here, $C_\sigma$ equals
$\binom{\w+1}{1}$ or $\binom{\w+2}{2}$
if $\sigma(1)$ is the first element of $J$ or not.
The coefficients $C_\sigma$ are used to correct for
multiplicities.

Consider the K\"unneth decomposition associated
to the first part of $P_\sigma$.
If $\sigma(1)$ is the first element of $J$, then
the terms with an odd class distributed to $\sigma(1)$ 
contribute to $V$ (and are normalized by the prefactor
$\binom{\w+1}{1}^{-1}$).
If the odd parts are distributed away from $\sigma(1)$,
then the resulting terms
are of the form \eqref{sominv}.
If $\sigma(1)$ is not the first element of $J$, then
all the K\"unneth distributions contribute to $V$
(and are normalized by the prefactor $\binom{\w+2}{2}^{-1}$).

We may then eliminate $V$ using the relations $R(0), \dots, R(|I|-1)$
with the induced descendent assignments and prefactor $M \tilde{\tau}(\omega)^\w$.
Because of the normalization and the removal of the
invariants \eqref{sominv}, the elimination analysis exactly follows
the proof of Lemma \ref{jjj}.
 
By degeneration, Lemma \ref{cakeshop}, and induction on $M$, we conclude 
the invariants
$$\lang M  \  \prod_{i\in I} \tilde{\tau}(\alpha) 
\ \tau_{m_1}(\beta) \prod_{1\neq j\in J}
\tilde{\tau}(\beta)\ | \ \eta \rang^E_d,$$
\begin{equation}\label{sominv3}
\lang M L \ \prod_{1\neq i\in I} \tilde{\tau}(\alpha) 
\ \tau_{m_1}(\beta)\prod_{1\neq\in J}
\tilde{\tau}(\beta) \rang^E_d,
\end{equation}
are inductively determined for any $m_1$.



We now analyze the assignment ${\mathbf l}_{f[r]l[s]}$
where $r+s > 1$. 
The
{\em special} elements are the first $r$ elements of $I$ and
the first $s$ elements of $J$.
Let $V$ be  the sum of the terms of
\begin{equation} \label{dd} 
\sum_\sigma C_\sigma^{-1}  
 (-1)^{\binom{|I|}{2}}  \text{sign}(\sigma) \ V_d(M,P_\sigma,{\mathbf l}_{f[r]l[s]})
\end{equation} 
with $|I|+|J|$ odd parts modulo the invariants 
determined by the analysis for the assignments ${\mathbf l}_{f[r']l[s']}$
for $r'+s'< r+s$.

In the definition of $V$, the summands are normalized with 
prefactors $C_\sigma^{-1}$ 
depending on the assignment function and $\sigma(1)$.
The first part $$\{1\} \cup W \cup \{\sigma(1)\}$$
of $P_\sigma$ contains either 2,1, or 0 special elements:
\begin{enumerate}
\item[$\bullet$]
if $P_\sigma$ contains 2 special elements, then  $C_\sigma= 1$,
\item[$\bullet$]
if $P_\sigma$ contains 1 special element, then  $C_\sigma =\binom{\w+1}{1}$,
\item[$\bullet$] if $P_\sigma$ contains 0 special elements, then
$C_\sigma =\binom{\w+2}{2}$.
\end{enumerate}

If $P_\sigma$ has special elements and
the distribution of odd classes
in the K\"unneth decomposition corresponding to $P_\sigma$ in a term 
of \eqref{dd} misses at least 1 special element,
then the term has {\em fewer} than $r+s$ special elements with odd classes.
Such terms are inductively determined by the analysis for
the assignments ${\mathbf l}_{f[r']l[s']}$ for $r'+s'  < r+s$.

We may then eliminate $V$ using the relations $R(0), \dots, R(|I|-1)$
with the induced descendent
assignments and prefactor $M \tilde{\tau}(\omega)^\w$.
Because of the normalization and the removal of the
invariants with fewer special elements, the elimination analysis exactly follows
the proof of Lemma \ref{jjj}.
 
Using degeneration, Lemma \ref{cakeshop}, and induction on $M$, 
the outcome for ${\mathbf l}_{f[r]l[s]}$
is a determination of all invariants 
$$\lang  M  \  
\prod_{i\leq r}\tau_{n_i}(\alpha) 
\prod_{r<i\in I} \tilde{\tau}(\alpha) \prod_{j\leq s} \tau_{m_j}(\beta)
\prod_{s<j\in J}
\tilde{\tau}(\beta)\ | \eta \rang^E_d,$$
$$
\lang  M L  \  \prod_{i\leq r}\tau_{n_i}(\alpha) \prod_{r<i\in I} 
\tilde{\tau}(\alpha) \prod_{j\leq s} \tau_{m_j}(\beta)
\prod_{s<j\in J}
\tilde{\tau}(\beta) \rang^E_d.$$
By induction on $r+s$, we find the invariant \eqref{yqpd} is
determined from the even theory by properties (i)-(iv).

The induction on $|I|$ is therefore established and 
the proof of Proposition \ref{bbbbb} is complete. 
\qed

\pagebreak

\section{Virasoro constraints for the full theory}

\label{three}

\subsection{Overview}
We complete the proofs of the main results of the paper.
Theorem \ref{gwh} was proven in \cite{OP2}. Theorems \ref{oddd} and
\ref{rv} are proven first. Theorem \ref{oddop} is then derived
as a corollary.

\subsection{Proof of Theorems \ref{oddd} and \ref{rv}}
We may define an {\em alternate}
relative theory of target curves by the following construction.
The alternate stationary sector is defined by the GW/H correspondence.
The descendents of the odd classes are added to the alternate theory
by the formula of Theorem \ref{oddd}.
The Virasoro constraints of Theorem \ref{rv} then define a unique extension
of the alternate theory including the descendents of the identity.
The proof of the existence and uniqueness of the Virasoro solution
here exactly follows the proof Proposition \ref{hhh}, the corresponding even result. 
The alternate theory of target curves is well-defined.

To prove Theorems \ref{oddd} and \ref{rv}, we must show the
alternate relative theory coincides with the relative Gromov-Witten theory.
Certainly, the two theories have equal stationary sectors by Theorem \ref{gwh}.
In fact, the two theories have equal even sectors since
we have proven the even relative Gromov-Witten satisfies the even
Virasoro constraints.

We now establish properties (i)-(iv), studied in Section \ref{threeodd}, hold for the
alternate relative theory of target curve.

\begin{enumerate}
\item[(i)] Algebraicity of the virtual class.

The balance of descendents of type $(1,0)$ and $(0,1)$ is the
only consequence of algebraicity used in Section \ref{threeodd}.
For odd classes in the
presence of descendents of $\omega$,
the alternate theory satisfies the balance property by the
formula of Theorem \ref{oddd}. Since 
the Virasoro constraints respect the
balance, the entire alternate theory satisfies the balance property.

\item[(ii)] Degeneration.

The GW/H correspondence is compatible with degeneration.
The formula of Theorem \ref{oddd} for the addition of the
odd classes is formally compatible with degeneration. 
The Virasoro constraints are also formally compatible with degeneration. 
Hence, the alternate theory satisfies the degeneration formula.

\item[(iii)] Monodromy invariance.

The stationary sector of the alternate theory is certainly
monodromy invariant. Since monodromy invariance preserves the
intersection form, the formula of Theorem \ref{oddd} for the
addition of the odd classes is monodromy invariant.
The monodromy invariance of the Virasoro solution is not 
immediate since a polarization of $H^*(X,\com)$ is
required for the definition of the Virasoro operators. However, 
an elementary argument by expansion in terms of the stationary theory
shows the elliptic monodromy relation,
$$R_d(N,{\mathbf n}, {\mathbf m},  \delta),$$
formally holds for the alternate theory.
Only these monodromy relations were used in Section \ref{threeodd}.

\item[(iv)] Elliptic vanishing relations.

An elementary argument by expansion in terms of the 
stationary theory shows the elliptic vanishing
relation,
$$V_d(M, P, {\mathbf l}),$$
formally holds for the alternate theory.
\end{enumerate}

In Section \ref{threeodd}, we proved the relative Gromov-Witten theory
of target curves is uniquely determined from the even theory by
properties (i)-(iv). Therefore, since the alternate theory coincides
with the relative Gromov-Witten theory on
the even sector and satisfies properties (i)-(iv), the alternate
theory equals the relative Gromov-Witten theory. 
\qed

\subsection{Proof of Theorem  \ref{oddop}}

Consider first the relative Gromov-Witten theory of target curves
{\em without} descendents of the identity. An explicit expansion
shows the differential equations 
\begin{eqnarray*}
D^i_{k} \ Z_d[\eta^1, \ldots, \eta^m] & =& 0, \\
\bar{D}_k^i \ Z_d[\eta^1, \ldots, \eta^m] & = & 0,
\end{eqnarray*}
when restricted to the zero locus of the ideal 
$$I=(t_0^0, t_1^0, t_2^0, \ldots),$$ are equivalent to
the formula of Theorem \ref{oddd}.

The differential equations for $D^i_k$ and $\bar{D}^i_k$
are proven to hold on the zero locus of
$I^r$ by induction on $r$ using the Virasoro constraints,
$$ L_n \ Z_d[\eta^1, \ldots, \eta^m] =0,$$
and the commutation relations,
\begin{align*}
           [L_n, D^i_k] & =  -(k+1) D^i_{n+k},\\
 [L_n, \bar{D}^i_k] & =  (n-k) \bar{D}^i_{n+k}.\\
\end{align*}
The Theorem is then deduced from completeness. \qed

\pagebreak

\vspace{+10 pt}
\noindent
Department of Mathematics \\
Princeton University \\
Princeton, NJ 08544\\
okounkov@math.princeton.edu \\

\vspace{+10 pt}
\noindent
Department of Mathematics\\
Princeton University\\
Princeton, NJ 08544\\
rahulp@math.princeton.edu
\end{document}